\documentclass[preprint]{elsarticle}

\usepackage{amssymb}
\usepackage{latexsym}

\usepackage{url}
\usepackage{xcolor}
\definecolor{newcolor}{rgb}{.8,.349,.1}

\usepackage[per-mode=symbol, exponent-product = \cdot]{siunitx}
\usepackage{bm}
\usepackage{amsmath}
\usepackage{cases}
\usepackage{amssymb}
\usepackage{physics}
\usepackage{float}
\usepackage{empheq}
\usepackage{hyperref}
\usepackage{cleveref}
\usepackage{todonotes}
\usepackage{comment}
\usepackage{graphicx} 
\usepackage{caption}
\usepackage{subcaption}
\usepackage{algorithm}
\usepackage{algorithmic}
\usepackage{mathrsfs}  
\usepackage{placeins}
\usepackage[margin=2cm]{geometry}
\usepackage{booktabs}
\usepackage{multirow}
\usepackage{pdflscape}
\usepackage{amssymb}
\usepackage{accents}
\usepackage{lineno}
\usepackage{soul}

\newcommand{\uale}{\bm u^\mathrm{ALE}}
\newcommand{\sigmak}{{\Sigma_\mathrm{k}}}
\newcommand{\struct}{\mathrm{S}}
\newcommand{\xhat}{{\widehat{\bm x}}}

\definecolor{ao(english)}{rgb}{0.0, 0.5, 0.0}
\newcommand{\rev}[1]{\textcolor{black}{#1}}

\renewcommand{\k}{{\mathrm{k}}}

\DeclareSIUnit\mmhg{mmHg}
\DeclareSIUnit\molar{M}

\begin{document}


\begin{frontmatter}

\title{An electromechanics-driven fluid dynamics model for the simulation of the whole human heart}%

\author[1,2]{Alberto {Zingaro}\corref{cor1}}
\cortext[cor1]{Corresponding author: email: alberto.zingaro@polimi.it}
\author[1]{Michele {Bucelli}}
\author[1]{Roberto {Piersanti}}
\author[1]{Francesco {Regazzoni}}
\author[1]{Luca {Dede'}}
\author[1,3]{Alfio {Quarteroni}}

\address[1]{MOX, Laboratory of Modeling and Scientific Computing, Dipartimento di Matematica, Politecnico di Milano, Piazza Leonardo da Vinci 32, 20133, Milano, Italy}
\address[2]{ELEM BioTech, Pier07, Calle Laietana, 26, 08003, Barcelona, Spain}
\address[3]{Institute of Mathematics, \'Ecole Polytechnique F\'ed\'erale de Lausanne, Station 8, Av. Piccard, CH-1015 Lausanne, Switzerland  (Professor Emeritus)}


\begin{abstract}
We introduce a multiphysics and  {geometric} multiscale computational {model, suitable to describe}  the hemodynamics of the whole human heart,  {driven by a four-chamber electromechanical model}.   {We {first} present a study on the calibration of the biophysically detailed RDQ20 \rev{active contraction model} (Regazzoni \textit{et al.}, 2020) that is able to reproduce the physiological range of hemodynamic biomarkers}.
{Then, w}e demonstrate that the ability of the force generation model to reproduce certain microscale mechanisms, such as the dependence of force on fiber shortening velocity, is crucial to  {capture}  the overall {physiological} mechanical and fluid dynamics  macroscale  behavior. This  {motivates} the need for using multiscale models with high biophysical fidelity, even when the outputs of interest are relative to the macroscale.  {We show that the {use} of a high-fidelity electromechanical model, combined with {a detailed} calibration process, allows us to achieve a remarkable biophysical fidelity {in terms of {both} mechanical and hemodynamic quantities}. Indeed, our electromechanical-driven CFD simulations {-- carried out on} an anatomically accurate geometry of the whole heart  {--}  provide results {that match the} cardiac physiology both qualitatively (in terms of flow patterns) and quantitatively (when comparing in silico results with biomarkers acquired in vivo).}
{Moreover}, we consider the pathological case of left bundle branch block, and we investigate the consequences that an electrical abnormality has on cardiac hemodynamics thanks to our multiphysics integrated model. The computational {model} {that} we propose can faithfully predict a delay and an increasing wall shear stress in the left ventricle {in the pathological condition}. The interaction of different physical processes in an integrated framework allows us to faithfully describe and model this pathology, {by} capturing and reproducing the intrinsic multiphysics nature of the human heart.
\end{abstract}


\end{frontmatter}


\section{Introduction}

The study of cardiac blood flows aims at enhancing the knowledge of heart physiology, assessing pathological conditions, and possibly improving clinical treatments and therapeutics. In the last decades, the role of mathematical models in the study of cardiac hemodynamics has increasingly gained relevance, for their non-invasiveness and flexibility with respect to geometries and flow conditions~\cite{chnafa2014image, this2020augmented, tagliabue2017complex, tagliabue2017fluid, Fumagalli2020, karabelas2018towards,this2020pipeline, masci2020proof, viola2020fluid}. Computational Fluid Dynamics (CFD) is largely employed to provide a detailed description of cardiac flows and to estimate hemodynamic indicators{, like e.g.,} the wall shear stress, that standard image-based techniques might not capture~\cite{chnafa2014image}.  {A biophysically detailed mathematical model of cardiac hemodynamics entails the complex interplay among different processes such as the interaction with cardiac electromechanics (EM), valve dynamics, transition-to-turbulence effects, and coupling with the surrounding circulation. Furthermore, while the literature about CFD models of the left ventricle~\cite{mangual2013comparative, zheng2012computational, seo2014effect, seo2013effect, tagliabue2017fluid, tagliabue2017complex},  {left atrium \cite{masci2017patient, masci2019impact, masci2020proof, zingaro2021hemodynamics, corti2022impact, vedula2015hemodynamics, koizumi2015numerical, dillon2018modeling, bosi2018computational, mazumder2022computational, zingaro2024comprehensive} } and left heart~\cite{viola2020fluid, viola2021effects, viola2022fsei, dede2021computational, zingaro2022geometric, chnafa2014image, zingaro2023comprehensive, bennati2023image, bennati2023turbulent, bennati2023image2} is {relevant}, the fluid dynamics simulation of the right heart is the subject of few works \cite{renzi2023accurate},  and they often focus on the sole right ventricle~\cite{wiputra2016fluid, collia2021comparative, mangual2012describing} {, neglecting the right atrium description}. Conversely, whole-heart fluid dynamics modeling is a much more challenging topic and has been addressed only recently in a few works  {\cite{mihalef2011patient, okada2019clinical, peirlinck2021precision, brenneisen2021sequential, karabelas2022global, viola2022gpu}.}}   {{As a matter of fact}, to the best of our knowledge, hemodynamics simulations of the whole heart are presented only in the following papers.} 
A four-chamber CFD model is proposed by the Siemens group in~\cite{mihalef2011patient}: the displacement of cardiac walls is obtained from patient-specific images and, since the left and right sides are not connected to the surrounding circulation, they performed simulations separately on the two parts. 
In~\cite{okada2019clinical}, the authors introduce a Fluid Structure Interaction (FSI) simulation of the whole heart based on the UT-Heart simulator developed at the University of Tokyo. 
Moreover, FSI simulations of the whole heart have been performed in the context of the Living Heart Project~\cite{peirlinck2021precision}.
A sequentially-coupled FSI model of the whole heart is devised in~\cite{brenneisen2021sequential}, showing that few iterations of the solver are needed to reach {convergence} in terms of mechanical indicators. 
Recently, a patient-specific whole-heart CFD model is introduced in~\cite{karabelas2022global}, consisting of the Navier-Stokes-Brinkman equations~\cite{fuchsberger2022incorporation} with prescribed boundary motion and, similarly to~\cite{mihalef2011patient}, simulations of the left and right parts are carried out separately. 
Finally, a GPU-accelerated fully-coupled electro-mechano-fluid computational model of the whole heart, based on the immersed boundary method, is presented in~\cite{viola2022gpu}.


{In this paper, we propose a  {four-chamber} electromechanical-driven fluid dynamics model, {which is characterized by} a remarkable biophysical fidelity and able to faithfully describe potential pathologies. Moreover, to  {improve the virtual representation of the heart physiology enabled by our model,} we present a \rev{careful manual} calibration of the \rev{active contraction model} (at the cellular level) to reproduce mechanical and hemodynamic biomarkers in the physiological range.
}

{The blood  {flow} in heart chambers is commonly modeled by Navier-Stokes equations for Newtonian fluids \cite{quarteroni2017cardiovascular}.} A crucial aspect in heart flows modeling is the treatment of boundary displacement, i.e. the way the deformation of cardiac walls is accounted for into the model. The boundary displacement can be the solution of a suitable mathematical model for the dynamics of the walls, fully coupled to the fluid dynamics model in an FSI framework, by imposing geometric, kinematic, and dynamic coupling conditions at the fluid-solid interface. The motion of the myocardium is in turn driven by muscular EM, resulting in a coupled electro-mechano-fluid problem (see e.g.~\cite{viola2020fluid, viola2021effects, viola2022fsei, choi2015new, santiago2018fully, santiago2018fluid, gerbi2018numerical, bucelli2022mathematical}). This approach, while being very comprehensive and physically motivated, entails a significant computational effort, due to the number of subsystems involved and to the non-linearity induced by the coupling. To mitigate this large computational cost, the boundary displacement can be prescribed as a datum, without any feedback from the fluid flow, in a CFD modeling framework. The displacement may be {prescribed by suitable} analytical laws~\cite{tagliabue2017fluid,dede2021computational, tagliabue2017complex, zingaro2021hemodynamics, domenichini2005three, baccani2002vortex, zheng2012computational, seo2013effect, corti2022impact, santiago2022design}, from patient-specific image-based reconstructions~\cite{Fumagalli2020, this2020pipeline, chnafa2014image, masci2019impact, masci2017patient, masci2020proof, bennati2023image, bennati2023image2, bennati2023turbulent, renzi2023accurate}, or {else} obtained from a previously performed EM simulation~\cite{augustin2016patient, karabelas2018towards, this2020augmented,  this2020augmented, zingaro2022modeling, zingaro2022geometric, zingaro2022mathematical}. The latter corresponds to a one-way coupled approach between EM and CFD, since only the kinematic coupling is enforced, without {foreseeing} any dynamic feedback from the fluid to the structure problem. This approach is also referred to as ``kinematic uncoupling''~\cite{this2020augmented, this2019image}.

A critical issue in 3D cardiovascular hemodynamics modeling is the prescription of boundary conditions at inlet and outlet sections, since boundary data are generally unavailable, but also because the circulatory system is a closed-loop network and the mathematical model should account for it. A possible approach is the so called \textit{geometric multiscale modeling}~\cite{quarteroni2016geometric}: the region of interest (in our case, the whole-heart) is described by a 3D model, while the remaining part of the circulation is addressed by means of lumped-parameter models, as 0D~\cite{blanco20103d, shi2006numerical, quarteroni2016geometric, milivsic2004analysis, kim2009coupling}, or 1D~\cite{quarteroni2016geometric, van2011pulse, formaggia2003one, formaggia2001coupling, santiago2018fluid}. The geometric multiscale modeling allows to account for the mutual interaction between the heart and the circulatory system, especially if the lumped parameter model provides a closed-loop description of the vascular network, as done in~\cite{zingaro2022geometric, fedele2022comprehensive, regazzoni2022cardiac, blanco20103d, augustin2021computationally, marcinno2022computational}.  

An additional key aspect in cardiac CFD simulations is the modeling of the cardiac valves. In principle, these can be treated by considering a structural model 
for the solid (leaflets of the valve and possibly its chordae tendinae) and a fluid dynamics model for the surrounding blood flow. This approach yields a coupled FSI model of the blood-valve system~\cite{ma2013image, kunzelman2007fluid, su2014numerical, gao2017coupled, de2003three, carmody2006approach, viola2022fsei, oks2022fluid, hiromi2021interface, spuhler20183d}, characterized by contact phenomena and fast dynamics. Thus, FSI valve models are commonly associated to a huge computational burden, to be added to the overall cost of the heart CFD simulation. To avoid this large computational cost, the effects of the valves in the blood can be surrogated by relying on reduced models for the valve dynamics~\cite{fedele2017patient, Fumagalli2020,fernandez2008numerical, astorino2012robust, this2019augmented, karabelas2022global, zingaro2022modeling}.

{Our computational model of the whole human heart encompasses} the main features {of the} cardiac hemodynamics: EM, cardiac valves, transition-to-turbulence effects, and {interplay} with the external circulation. Indeed, our fluid model is driven by the four-chamber EM model  recently proposed in~\cite{fedele2022comprehensive}. The multiphysics model is fully coupled to the external circulation described by a lumped-parameter model, extending the computational framework we introduced in \cite{zingaro2022geometric} to the case of four-chamber CFD simulations. We carry out numerical simulations on an anatomically accurate geometry of the heart, obtaining results that are quantitatively in agreement with data from the medical literature and qualitative faithfully in terms of blood flow patterns. In this respect, we analyze the role played by the highly biophysically detailed RDQ20 \rev{active contraction model}~\cite{regazzoni2020biophysically} {on relevant} hemodynamic quantities.  {Indeed, since the electromechanical displacement drives the deformation of cardiac chambers for the CFD simulation, its calibration is fundamental towards faithfully reproducing heart physiology. As shown in \cite{fedele2022comprehensive}, the parameters of the \rev{active contraction model} play a significant role in determining the flow rates across cardiac valves, which have a dramatic impact on the CFD simulation, both in terms of macroscopic indicators and overall flow distribution. Therefore, we present a detailed \rev{manual} calibration of the active force generation model, validating our results  {on} several macroscopic heartbeat indicators such as stroke volumes, ejection fractions, peak flowrates and, consequently, also in terms of blood velocities computed in the CFD simulation. Our  {sensitivity analysis} highlights  {that} microscopic features of the RDQ20 model have a large impact on the macroscopic characteristics of the heartbeat.}
Then,  {we show that our detailed computational model can  {improve the understanding of the impact of  the Left Bundle Branch Block (LBBB) pathology {\cite{tan2020left}} on CFD biomarkers}, by capturing the effects that an electrophysiological abnormality has on different physical processes behind the heart activity, {henceforth} allowing to capture the intrinsic multiphysics nature of the cardiac function. } { This paper represents one the few examples in the literature on the modeling and simulation of the whole heart hemodynamics. Moreover, to the best of our knowledge, this is the first work in which the 3D whole heart fluid dynamics model is also coupled to a lumped-parameter closed-loop circulation model. }

This paper is organized as follows: in \Cref{sec:mathematical-models}, we introduce the mathematical models employed for the EM, fluid dynamics and circulation problems. \Cref{sec:numerical-methods} is devoted to the description of the numerical methods for each subproblem and to the strategies used for their coupling. In \Cref{sec:numerical-results}, we present numerical results on a realistic whole-heart geometry, in {both} physiological and pathological conditions. Finally, limitations and conclusions follow in \Cref{sec:limitations} and \Cref{sec:conclusions}, respectively.

\section{Mathematical model}
\label{sec:mathematical-models}
In this section, we introduce the mathematical model. Specifically, the EM model is briefly described in \Cref{sec:model_em} and the whole-heart fluid domain is defined in \Cref{sec:wh-fluid-domain}. We present {our approach suitable} to deal with the domain deformation in \Cref{sec:model_geo}, the fluid dynamics model in \Cref{sec:modeling_fluid},  {and its coupling with the external circulation in} \Cref{sec:coupling_circulation}.

\subsection{The whole-heart electromechanical model}
\label{sec:model_em}

\begin{figure}[t]
	\centering
	\includegraphics[trim={1 1 1 1},clip,width=\textwidth]{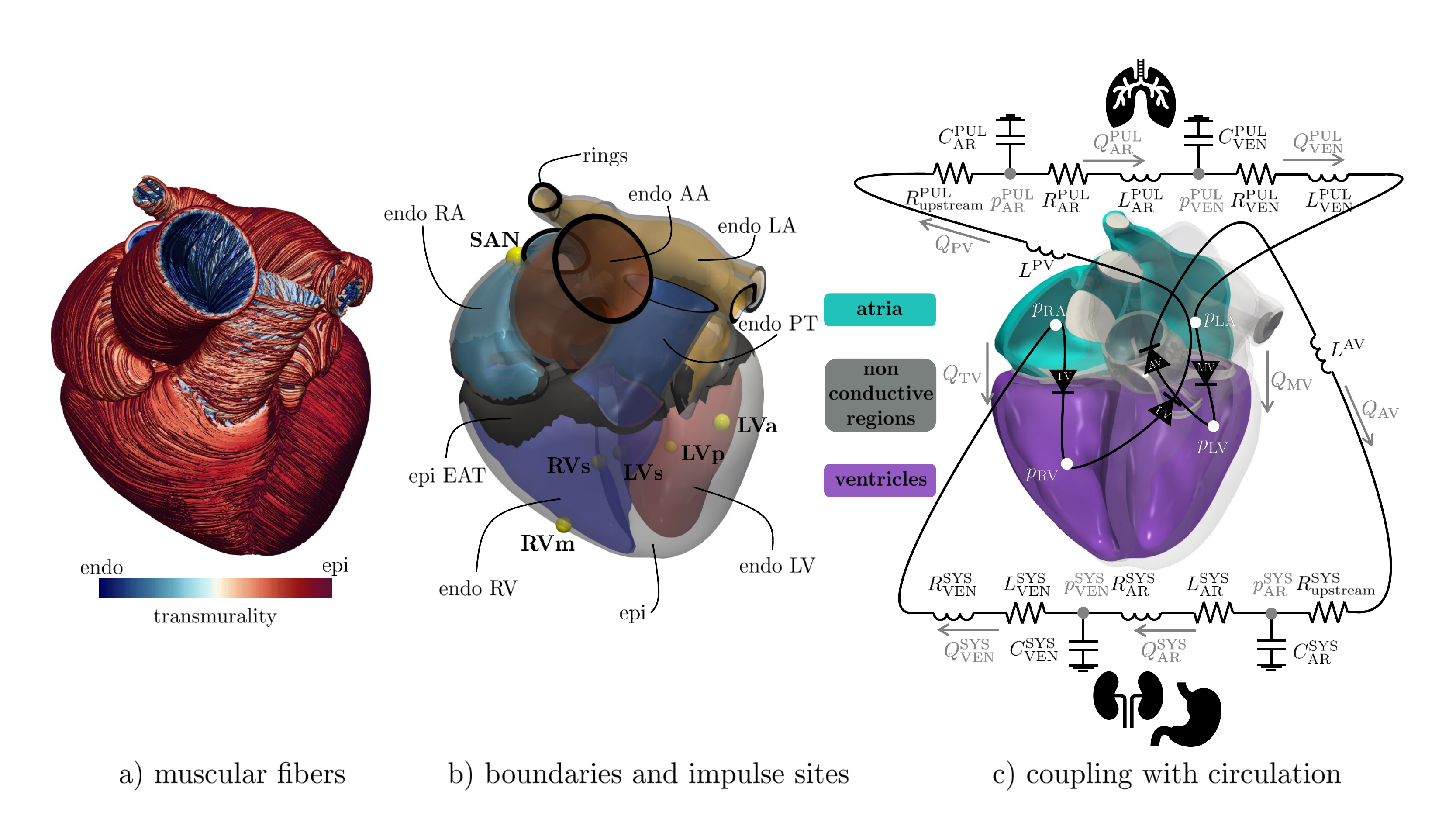}
	\caption{The whole-heart EM model: a) cardiac muscular fibers, (b) boundaries and impulse sites (yellow spheres with bold labels), (c) coupling with circulation and we highlight the three main regions of the EM model (atria, ventricles and non-conductive regions).}
	\label{fig:em-circulation-epicardium-fibers}
\end{figure}

The whole-heart EM is based on a comprehensive and biophysically detailed computational model that we recently presented in~~\cite{fedele2022comprehensive}. More precisely, we consider a 3D description of cardiac EM in all the four-chambers and a 0D representation of the complete circulatory system, including the cardiac blood hemodynamics~\cite{regazzoni2022cardiac,piersanti2022biv}, as we display in  \Cref{fig:em-circulation-epicardium-fibers}c. 

The whole-heart EM model includes a detailed myocardial fiber architecture built upon a total-heart Laplace-Dirichlet Rule-Based Method~\cite{piersanti2021phd}, which couples together different methods for the atria~\cite{piersanti2021fibers} and the ventricles~\cite{doste2019rule}, to properly reproduce the characteristic features of the cardiac fiber bundles in all the four-chambers~\cite{piersanti2021fibers}, see \Cref{fig:em-circulation-epicardium-fibers}a.

Cardiac electrophysiology is described by means of the monodomain equation equipped with no-flux Neumann boundary conditions {\cite{franzone2014mathematical}} and endowed with the following human ionic models: ten Tusscher-Panfilov for the ventricles~\cite{TTP06} and Courtemanche-Ramirez-Nattel for the atria~\cite{CRN}. Furthermore, the arterial vessels and the atrio-ventricular basal plane are assumed to be non-conductive regions, {whence} electrically isolating the atria from the ventricles~\cite{fedele2022comprehensive}. Finally, the cardiac conduction system is substituted by a series of spherical electrical impulses{, originating} from the sino-atrial node (SAN) and ending into the left and right ventricular endocardia which, combined with a fast endocardial layer, surrogates the effect of the Purkinje network~~\cite{lee2019fastendo,fedele2022comprehensive}, as we show in  \Cref{fig:em-circulation-epicardium-fibers}b.

The sarcomere mechanical activation is based on the biophysically detailed RDQ20 active contraction model~\cite{regazzoni2020biophysically}, properly calibrated for both atria~\cite{mazhar2021} and ventricles~\cite{regazzoni2020activation}. The RDQ20 is able to represent in detail the sophisticated microscopic active force generation mechanisms, taking place at the scale of sarcomeres~\cite{regazzoni2020activation}. Moreover, to differentiate the active tension in left and right ventricles, we consider a spatially heterogeneous active tension~\cite{piersanti2022biv} \rev{, i.e. defined in relation to a normalized inter-ventricular distance. Notably, we set the active tension in the right ventricle at 70\% of that in the left ventricle, whereas we keep the same contractility between left and right atrium.}

The myocardial tissue mechanics is described by the momentum balance equation under the hyperelasticity assumption~\cite{ogden1997}. We employ, for the active part, an orthotropic active stress formulation~\cite{piersanti2022biv}, which surrogates the contraction caused by dispersed myofibers~\cite{guan2020dispers}, and, for the passive behavior, specific mechanical constitutive laws and model parameters for the different cardiac region: the Usyk constitutive law for both the atria and the ventricles~\cite{usyk2002} and a Neo-Hookean strain energy density function for the atrio-ventricular basal plane and the vessels~\cite{ogden1997}. Finally, a nearly incompressible formulation is enforced with a penalty method~\cite{fedele2022comprehensive}. Concerning the mechanical boundary conditions, we consider: i) generalized Robin boundary conditions on the epicardium, surrogating both the presence of the pericardium and also the epicardial adipose tissue, crucial for reproducing the correct downward and upward movement of the atrio-ventricular basal plane~\cite{fedele2022comprehensive}; ii) normal stress boundary conditions on the four-chamber endocardia and vessel endothelia to account for the pressures exerted by the blood, where the endocardium and endothelium fluid pressures are given by the coupling between the mechanical and the circulation problems~\cite{regazzoni2022cardiac,piersanti2022biv,fedele2022comprehensive}; iii) homogeneous Dirichlet boundary condition on all the artificial rings where we cut the computational domain, since the arteries and atrial veins can be considered fixed here~\cite{fedele2022comprehensive}, as displayed in \Cref{fig:em-circulation-epicardium-fibers}b. \rev{When setting boundary conditions on the epicardial boundary, we opted for reduced values of normal stiffness and normal damping on the epicardial adipose tissue (EAT, see \Cref{fig:em-circulation-epicardium-fibers}b). This choice allows the ventricular base to be more free to move, as well as the appendages of the left and the right atria.}

The whole-heart 3D EM model is fully coupled with a 0D closed-loop lumped parameters model for the blood hemodynamics through the entire cardiovascular network. Systemic and pulmonary circulations are modeled {using} resistance-inductance-capacitance circuits (both for the arterial and venous part) and non-ideal diodes stand for the heart valves~~\cite{regazzoni2022cardiac}. In \Cref{fig:em-circulation-epicardium-fibers}c we give a graphical representation of the 3D-0D model. The coupling between the 0D and 3D EM models is achieved by introducing volume-consistency coupling conditions, where the pressures of all the four-chambers act as Lagrange multipliers associated with the introduced volume constraints~~\cite{regazzoni2022cardiac,piersanti2022biv}. 

The most relevant feedbacks, representing the interactions among electric signal propagation, the cardiac tissue deformation and contraction, and the circulatory system{,} are modelled inside the whole-heart EM model~\cite{fedele2022comprehensive}. These include e.g. the mechano-electric feedback~\cite{salvador2022mef} (between electrophysiology and mechanics) and the fibers-stretch and fibers-stretch-rate feedbacks~\cite{regazzoni2019reviewXB} (between mechanics and the \rev{active contraction model}).  

For {the full set of equations} of the whole-heart EM model, we refer to~\cite{fedele2022comprehensive,piersanti2022biv,regazzoni2022cardiac}. 

{In this paper,} the whole-heart EM model serves as unidirectional input for the fluid dynamics problem as we better detail in \Cref{sec:model_geo} and \Cref{sec:modeling_fluid}.

\subsection {The whole-heart fluid domain}
\label{sec:wh-fluid-domain}
\begin{figure}[!t]
	\centering
	\includegraphics[trim={1 3cm 1 1 },clip,width=\textwidth]{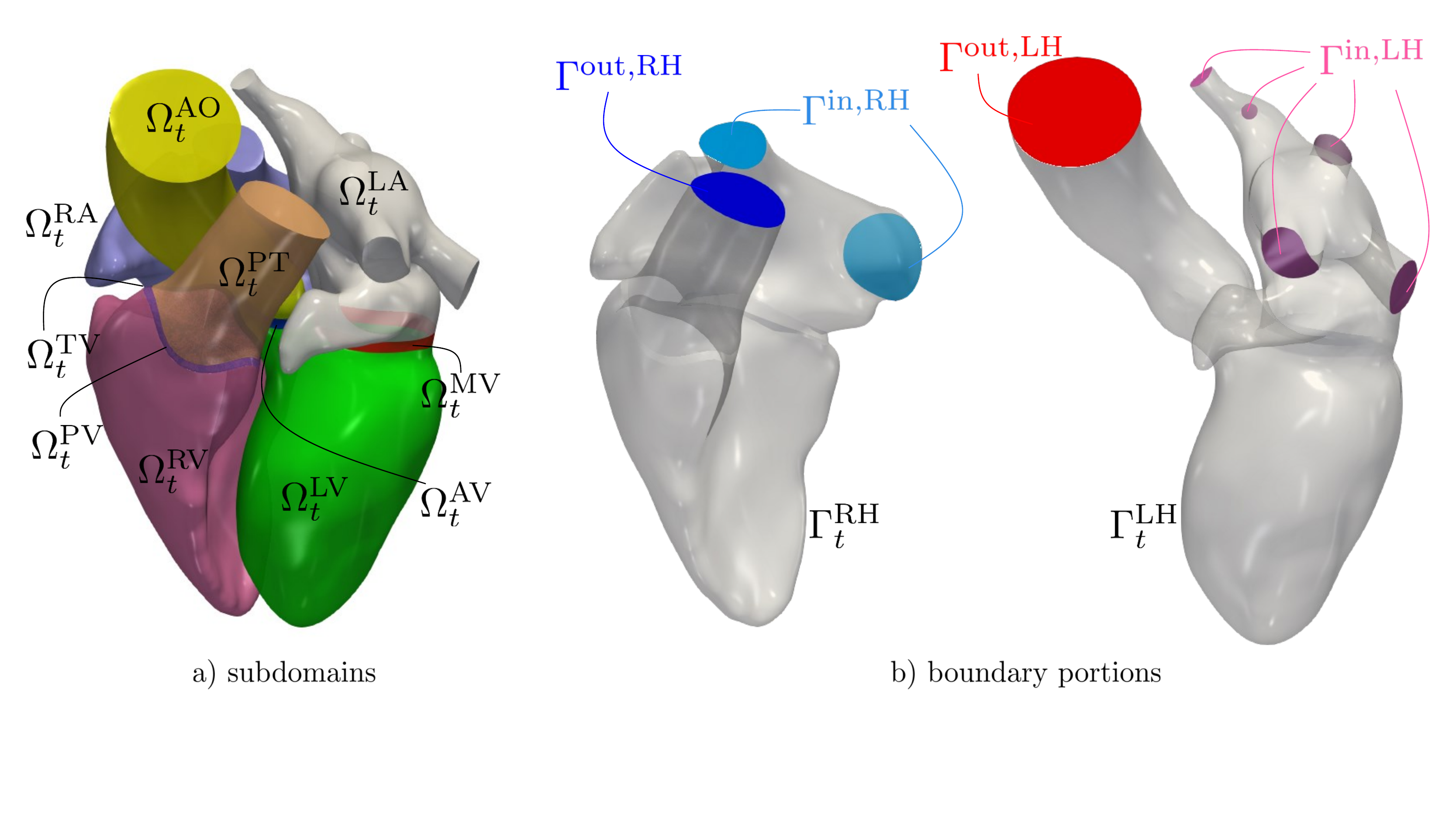}
	\caption{The whole-heart fluid domain: (a) subdomains composing the whole heart; (b) boundary portions (the left and right part are separated for visualization purposes).}
	\label{fig:domain_tags}
\end{figure}
Let $\Omega_t$ be the fluid domain at time $t>0$ bounded with a sufficiently regular boundary $\Gamma_t \equiv \partial \Omega_t $ and let $(0, T)$ be the temporal domain, with $T$ the final time. From a fluid dynamics view point, the whole-heart fluid domain $\Omega_t$ is topologically disjoint and split into left heart (LH, $\Omega_t^\mathrm{LH}$) and right heart (RH, $\Omega_t^\mathrm{RH}$), as we show	 in \Cref{fig:domain_tags}:  $
\Omega_t = \Omega_t^\text{RH} \cup \Omega_t^\text{LH} $,  with $\overline{\Omega}_t^\text{RH} \cap \overline{\Omega}_t^\text{LH} = \varnothing$. Specifically,
\begin{align*}
	{\overline{\Omega}_t^\text{RH}} & {=  \overline{\Omega}_t^{\text{RA}} \cup \overline{\Omega}_t^{\text{TV}}\cup \overline{\Omega}_t^{\text{RV}} \cup \overline{\Omega}_t^{\text{PV}} \cup \overline{\Omega}_t^{\text{PT}}, }
	\\
	{\overline{\Omega}_t^\text{LH}} &  {=  \overline{\Omega}_t^{\text{LA}} \cup \overline{\Omega}_t^{\text{MV}}  \cup \overline{\Omega}_t^{\text{LV}} \cup \overline{\Omega}_t^{\text{AV}} \cup \overline{\Omega}_t^{\text{AO}}},
\end{align*}
where  $\Omega_t^{\text{RA}}$, $\Omega_t^{\text{RV}}$, $\Omega_t^{\text{PT}}$ are the right atrium (RA), right ventricle (RV) and pulmonary trunk (PT) subdomains, and $\Omega_t^{\text{LA}}$, $\Omega_t^{\text{LV}}$, $\Omega_t^{\text{AO}}$ the left atrium (LA), left ventricle (LV) and aorta (AO) subdomains. Moreover, $\Omega_t^{\text{TV}}$,  $\Omega_t^{\text{PV}}$, $\Omega_t^{\text{MV}}$, $\Omega_t^{\text{AV}}$ are the subdomains representing the rings of the tricuspid valve (TV), pulmonary valve (PV), mitral valve (MV) and aortic valve (AV). Analogously, we partition the boundary of the whole-heart domain as $ \Gamma_t = \Gamma_t^\text{RH} \cup \Gamma_t^\text{LH} $, with {$\Gamma_t^\text{RH} \cap \Gamma_t^\text{LH} = \varnothing$}. In particular, as displayed in \Cref{fig:domain_tags}b,
\begin{equation*}
	\Gamma_t^\text{RH} = \Gamma^\text{out,RH} \cup \Gamma^\text{in,RH} \cup \Gamma_t^{\text{w,RH}},
\end{equation*} 
with $\Gamma^\text{in,RH}$ the inlet sections of the superior and inferior venae cavae, $\Gamma^{\text{out,RH}}$ the outlet section of the pulmonary trunk, and $\Gamma_t^{\text{w,RH}}$ the endocardium of the RH. In an analogous fashion, on the left part:
\begin{equation*}
	\Gamma_t^\text{LH} = \Gamma^\text{out,LH} \cup \Gamma^\text{in,LH} \cup  \Gamma_t^{\text{w,LH}},
\end{equation*} 
with $\Gamma^\text{in,LH}$ the five inlet sections of the four pulmonary veins, $\Gamma^\text{out,LH}$ the outlet section of the aorta, and $\Gamma_t^{\text{w,LH}}$ the LH endocardium. 

\subsection{The {fluid domain displacement} problem}
\label{sec:model_geo}
To represent the deformation of the domain over time, we introduce a fixed reference configuration $\hat\Omega \subset \mathbb{R}^3$, such that the domain in current configuration $\Omega_t$ is defined at any $t \in (0, T)$ as
\begin{equation*}
	\Omega_t  = \{ \bm x \in \mathbb R^3: \; \bm x = \xhat + \bm d(\xhat, t), \; \xhat  \in \widehat \Omega \},
\end{equation*}
where $\bm d: \widehat{\Omega} \times (0, T)$ is the displacement of the domain{, and is obtained} it by solving the following harmonic extension problem: 
\begin{subnumcases}{\label{eq:lifting}}
	- \div (s \grad \bm d) = \mathbf 0 &  in $\widehat \Omega \times (0, T),$ \\
	\bm d = \bm d^{\partial \Omega}(\mathbf x, t) &  on $ {\partial \widehat \Omega} \times (0, T)$.
\end{subnumcases}
In \Cref{eq:lifting}, $\bm d^{\partial \Omega}: {\partial \widehat \Omega} \times (0, T)$ is the boundary displacement, computed by restricting the solution of the EM simulation to the endocardium and the endothelium. Furthermore, $s: {\widehat \Omega} \times (0, T) \to \mathbb R$ is a space-dependent scalar field introduced to avoid distortion of mesh elements.  Specifically, we use the boundary-based stiffening approach proposed in~\cite{jasak2006automatic}. We denote the {fluid domain displacement} problem \Cref{eq:lifting}  with the abridged notation
\begin{equation*}
	{\mathscr D}(\bm d, \bm d^{\partial \Omega}) = 0.
\end{equation*}
The EM simulation is solved with a significantly larger timestep than the CFD one. Therefore, the boundary displacement $\bm d^{\partial\Omega}$ is only available for some times $t_k$, $k = 0, 1, \dots, N_\text{EM}$, although the domain displacement $\bm d$ is needed with a finer temporal resolution. Thus, problem \eqref{eq:lifting} is solved for all times $t_k$, and then we construct a displacement field $\tilde{\bm d}(t) : \widehat\Omega \times (0, T) \to \mathbb{R}^3$ using smoothing splines approximation in time ~\cite{de1978practical}.

We compute the domain velocity by deriving the displacement in time as
\begin{equation}
	\uale = \pdv{\widetilde{\bm d}}{t} \; \; \mathrm{in} \; \Omega_t \times (0, T). 
	\label{eq:velocity_ALE}
\end{equation}

\subsection{The Navier-Stokes equations in ALE framework with {the} RIIS {model} of {the} valves}
\label{sec:modeling_fluid}
We {model} the blood in the cardiac cavities as an incompressible, viscous and Newtonian fluid characterized by constant density $\rho$ and constant dynamic viscosity $\mu$. {We therefore use} the time-dependent incompressible Navier-Stokes equations expressed in an Arbitrary Lagrangian Eulerian (ALE) framework to account for the moving domain. We denote by $\bm u:\Omega_t \times (0, T) \to \mathbb R^3 $ and $  p:\Omega_t \times (0, T) \to \mathbb R$ the fluid velocity and pressure, respectively. Let $\sigma$ be the Cauchy stress tensor, defined for incompressible, Newtonian and viscous fluids as $\sigma (\bm u, p) = -p I + 2 \mu \epsilon(\bm u)$, with $\epsilon(\bm u) = \frac{1}{2}\left (\grad \bm u + (\grad \bm u)^T \right )$ the strain-rate tensor. 

We model the effects of cardiac valves in the fluid by means of the Resistive Immersed Implicit Surface (RIIS) method~\cite{fedele2017patient}. We consider four immersed surfaces $\sigmak$, with $\k \in \mathcal{I}_{\mathrm v} = \left \{ \text{MV, AV, TV, PV} \right \}$ the set of valves. Each valve is characterized by a resistance coefficient $R_\k$ and a parameter $\varepsilon_\k$ representing the half thickness of the valve leaflets. The immersed surface is implicitly described by a signed distance function $\varphi_\k: \Omega_t \times (0, T) \to \mathbb R$. With the RIIS method, we introduce the following penalty term to the momentum balance of the Navier-Stokes equations (expressed in ALE form):
\begin{equation*}
	\bm{\mathcal R}(\bm u, \uale) = \sum_{\k \in \mathcal{I}_{\mathrm v}} \frac{R_\mathrm{k}}{\varepsilon_\mathrm{k}} \delta_{\sigmak, \varepsilon_\mathrm{k}}(\varphi_\mathrm{k})\left(\bm u - \uale - \bm u_{\sigmak}\right ).
	\label{eq:riis}
\end{equation*}
$\bm {\mathcal R}$ penalizes the mismatch between the relative velocity $\bm u - \uale$ and the velocity of the valves' leaflets $\bm u_{\sigmak}$, weakly imposing a kinematic coupling condition only. The resistive term {is only acting on a tiny support} around $\sigmak${, thanks to}  {a} smoothed Dirac delta function {acting as multiplicative factor}.  {We refer to \cite{fedele2017patient} for its definition.}

By defining with $\frac{\widehat \partial \bm u }{\partial t} =  \pdv{\bm u}{t} + (\uale \cdot \nabla ) \bm u  $ the ALE derivative, the 3D fluid dynamics model of the whole heart reads: 
\begin{subnumcases}{\label{eq:ns_ale_riis}}
	\rho \frac{\widehat \partial \bm u}{\partial t} + ((\bm u - \uale)\cdot \grad )\bm u+ \div \sigma (\bm u, p) 
	+ \bm {\mathcal R}(\bm u, \uale)  = \bm 0
	& in  $\Omega_t \times (0, T)$,
	\\
	\div \bm u  = 0 &  in $ \Omega_t \times (0, T)$,
	\\
	\sigma (\bm u, p) \bm n  =  -p^\mathrm{in}_\mathrm{RA} \bm n & on $ \Gamma^\mathrm{in,RH} \times(0, T)$,
	\\
	\sigma (\bm u, p) \bm n  = -p^\mathrm{PUL}_\mathrm{AR} \bm n & on $ \Gamma^\mathrm{out,RH} \times(0, T),$
	\\
	\sigma (\bm u, p) \bm n  =  - p^\mathrm{in}_\mathrm{LA} \bm n & on $ \Gamma^\mathrm{in,LH} \times(0, T),$
	\\
	\sigma (\bm u, p) \bm n  = -p^\mathrm{SYS}_\mathrm{AR}  \bm n & on  $  \Gamma^\mathrm{out,LH} \times(0, T),$
	\\
	\bm u  =  \uale  &  on  $\Gamma_t^\text{w,RH} \cup \Gamma_t^\text{w,LH} \times (0, T). $
	\\
	\bm u = \bm u_0 & in $\Omega_0 \times  \{0\}$,
\end{subnumcases}
where $p^\mathrm{in}_\mathrm{RA}, \, p^\mathrm{PUL}_\mathrm{AR},\,  p^\mathrm{in}_\mathrm{LA}, \, p^\mathrm{SYS}_\mathrm{AR}  $ are the pressures arising from the coupling with the circulation model, as we detail in \Cref{sec:coupling_circulation}, and $\bm u_0$ is the initial velocity. We denote the 3D fluid dynamics model of the whole heart in \Cref{eq:ns_ale_riis} by
\begin{equation*}
	\mathscr F(\bm u, p, \uale) = 0.
\end{equation*}

\subsection{Coupling with circulation}
\label{sec:coupling_circulation}
\begin{figure}[t]
	\centering
	\includegraphics[trim={10mm 1cm 2 2},clip,width=\textwidth]{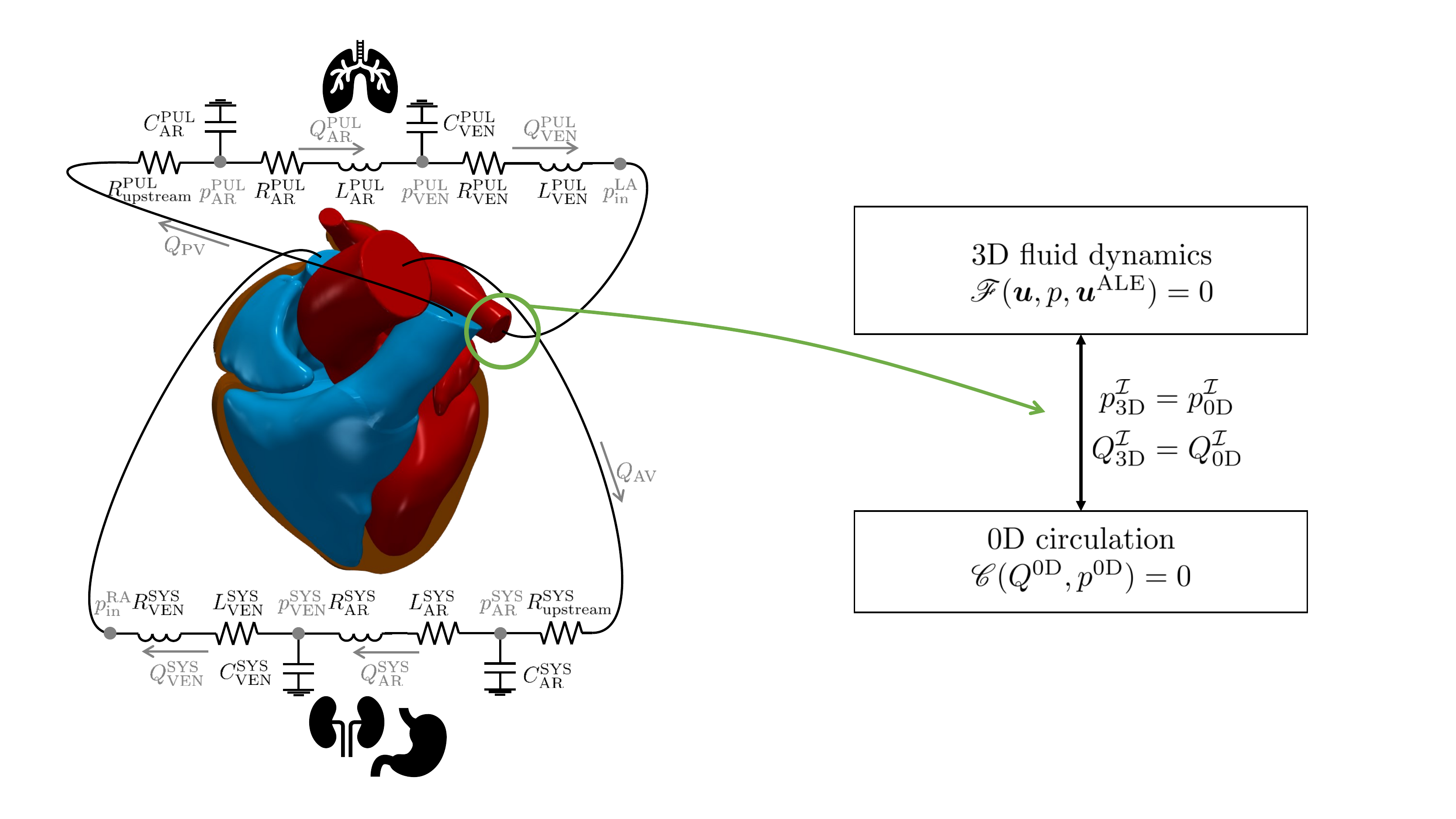}
	\caption{The 3D-0D fluid dynamics model of the whole-heart coupled to the surrounding circulation.}
	\label{fig:whole_heart_3d_0d}
\end{figure}
To account for the {interplay} between the heart's fluid dynamics and the hemodynamics of the surrounding circulation, we couple the cardiac CFD model to a 0D closed-loop model of the whole circulation proposed in~\cite{regazzoni2022cardiac}.  Specifically, we extend to the whole heart the coupling strategy that we devised in~\cite{zingaro2022geometric} for the case of the sole left heart. 
We consider a reduced (open) version of the original circulation model in which we remove the equations for all the variables that are already described by the 3D counterpart. The equations of the open system are reported in \ref{appendix:0dmodel}. By denoting with $\bm Q^\mathrm{0D}$ and $\bm p^\mathrm{0D}$ the vectors containing flowrates and pressures of the 0D model, respectively, we refer to the open system with the notation 
\begin{equation*}
	\mathscr C (\bm Q^\mathrm{0D}, \bm p^\mathrm{0D}) =  0.
\end{equation*}
The coupling between the 3D {CFD model} and the open 0D model consists in the enforcement of the continuity of pressures and flowrates on the artificially chopped boundaries $\Gamma^\mathcal{I}$ of the fluid domain, yielding the following  conditions:
\begin{subnumcases}{\label{eq:interfaces_conditions}}
	p^\mathcal{I}_\mathrm{3D} = p^\mathcal{I}_\mathrm{0D} &  $\text{ on } \Gamma^\mathcal{I} \times(0, T),$
	\label{eq:interfaces_conditions_dynamic}
	\\
	Q^\mathcal{I}_\mathrm{3D} = Q^\mathcal{I}_\mathrm{0D}  &  $\text{ on } \Gamma^\mathcal{I} \times(0, T),$
	\label{eq:interfaces_conditions_kinematic}
\end{subnumcases}
which express dynamic and kinematic coupling, respectively, with\footnote{We define the sign of the flowrate in accordance with the outward unit normal $\bm n$. Thus, an inlet flowrate (entering velocity) will be, by definition, negative.} 
\begin{equation}
	p^\mathcal{I}_\mathrm{3D} =  \frac{1}{|\Gamma^\mathcal{I}|}\int_{\Gamma^{\mathcal I}} p, \qquad
	Q^\mathcal{I}_\mathrm{3D} = \int_{\Gamma^{\mathcal I}} (\bm u  - \uale) \cdot \bm n.
	\label{eq:definition-p3D-Q3D}
\end{equation}
Considering the whole-heart fluid domain (see \Cref{fig:domain_tags}b), the interface boundaries are $\Gamma^\mathcal{I} = \Gamma^\mathrm{in, RH} \cup \Gamma^\mathrm{out, RH} \cup \Gamma^\mathrm{in, LH} \cup \Gamma^\mathrm{out, LH}$. The 0D pressures ($p^\mathcal{I}_\mathrm{0D}$) and flowrates ($p^\mathcal{I}_\mathrm{0D}$) are~\cite{bucelli2022mathematical, janela2010absorbing}:
\begin{equation*}
	\begin{aligned}
		& p^\text{in,RH}_\mathrm{0D} = p_\mathrm{RA}^\mathrm{in}, 
		& & Q^\text{in,RH}_\mathrm{0D} = - Q^\mathrm{SYS}_\mathrm{VEN}, 
		\\
		& p^\text{out,RH}_\mathrm{0D} = p^\mathrm{PUL}_\mathrm{AR} + R^\mathrm{PUL}_\mathrm{upstream}Q_\mathrm{PV}, 
		& & Q^\text{out,RH}_\mathrm{0D} = Q_\mathrm{PV}, 
		\\
		& p^\text{in,LH}_\mathrm{0D} = p^\mathrm{in}_\mathrm{LA}, 
		& & Q^\text{in,LH}_\mathrm{0D} =- Q^\mathrm{PUL}_\mathrm{VEN} ,
		\\
		& p^\text{out,LH}_\mathrm{0D} = p^\mathrm{SYS}_\mathrm{AR} + R^\mathrm{SYS}_\mathrm{upstream}Q_\mathrm{AV}, 
		& & Q^\text{out,LH}_\mathrm{0D} = Q_\mathrm{AV}.
	\end{aligned}
\end{equation*}
From the point of view of the 3D CFD model, the conditions expressed by \Cref{eq:definition-p3D-Q3D} are defective, since they prescribe the average pressure and the total flow rate over the entire section $\Gamma^\mathcal{I}$, rather than pointwise stress and velocity distributions~\cite{quarteroni2016geometric}. 
We choose to complete the pressure condition as
\begin{equation}
	\sigma(\bm u, p) \bm n = - p^\mathcal{I} \bm n, \text{ on } \Gamma^{\mathcal I} \times (0, T).  
	\label{eq:bc-neumann}
\end{equation} 
The flowrate condition, conversely, is left in its defective form, since it is sufficient for the algorithm we use for the 3D-0D coupling (see \Cref{sec:numerical-methods}).


\section{Numerical methods}
\label{sec:numerical-methods}

\begin{algorithm}[t]
	\caption{Numerical scheme for the EM-driven CFD simulation of the whole heart}
	\label{algo:coupling}
	\begin{algorithmic}
		\STATE{Solve whole-heart EM model}
		\STATE {Pick solution in the last heartbeat: $\to \bm d^\mathrm{EM}_{i}$, for $i = 0, \dots, N$}
		\STATE{Restrict $\bm d^\mathrm{EM}_{i}$ on $\partial \Omega_{i}^{\struct, \mathrm{endo}}$: $
			\bm d^{\partial \Omega}_{i}, \text{ for } i = 0, \dots, N$
		}   
		\STATE{Solve {fluid domain displacement} problem: ${\mathscr D} (\bm d_{i}, \bm d^{\partial \Omega}_{i}) = 0, \text{ for } i = 0, \dots, N$}
		\STATE{Build approximant: $\widetilde{\bm d}(t)$}
		\STATE{Initialization: $n= 0, \, \bm{u}_0 = \bm 0$}
		\WHILE{$n < N_t$}
		\STATE{
			{Update valves leaflet position}
			\\
			Compute ALE velocity: $\uale_{n+1} = \frac{\widetilde{\bm d}_{n+1} - \widetilde{\bm d}_{n}}{\Delta t}$.
			\\
			{Solve circulation}: $ \mathscr C (\bm Q^\mathrm{0D}_{n+1}, \bm p^\mathrm{0D}_{n+1}) = 0$ with data $ Q^\mathcal{I}_{\mathrm{0D}, n}$
			\\
			{Compute interface data (0D $\to$ 3D)}: $p^\mathcal{I}_{\mathrm{0D}, n+1} \to p^\mathcal{I}_{\mathrm{3D}, n+1}$
			\\
			{Solve fluid dynamics:} $ \mathscr F(\bm u_{n+1}, p_{n+1}, \uale_{n+1}) = 0$, with Neumann data $ p^\mathcal{I}_{\mathrm{3D}, n+1}$
			\\
			{Compute interface data  (3D $\to$ 0D)}: $Q^\mathcal{I}_{\mathrm{3D}, n+1} \to Q^\mathcal{I}_{\mathrm{0D}, n+1}$
			\\
			$n \leftarrow n+1.$ 
		}
		\ENDWHILE
	\end{algorithmic}
\end{algorithm}

\begin{figure}
	\centering
	\includegraphics[trim={1 6cm 1 1},clip,width=\textwidth]{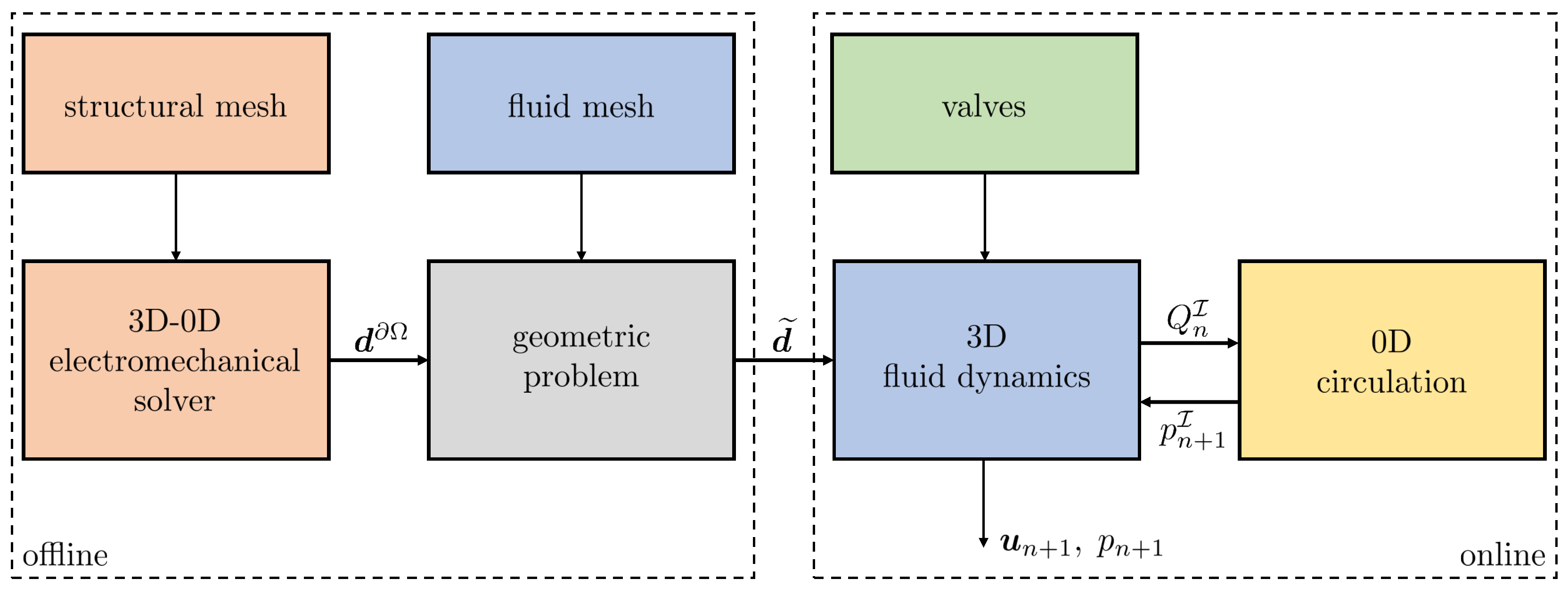}
	\caption{Graphical representation of the overall algorithm to simulate the whole-heart hemodynamics driven by the EM model.}
	\label{fig:algorithm}
\end{figure}

In {this section}, we describe the numerical methods we use to solve our multiphysics and multiscale system. The overall algorithm is presented in \Cref{algo:coupling} and graphically represented in \Cref{fig:algorithm}. We can subdivide the overall procedure in a {preliminary phase, in which we solve the EM simulation~\cite{fedele2022comprehensive} and the fluid domain displacement problem to lift the boundary displacement to the fluid bulk domain, followed by the coupled 3D-0D CFD simulation}.

\paragraph{Numerical methods for the EM model}
For the numerical approximation of the whole-heart EM model we employ the efficient Segregated-Staggered scheme~\cite{regazzoni2022cardiac,piersanti2022biv,fedele2022comprehensive}.
In this numerical scheme, the different cardiac physical models, contributing to both the 3D EM and the 0D blood circulation, are sequentially solved in a segregated manner, using different resolutions in space and time to properly account for the heterogeneous space and time scales characterizing different {physical processes}~\cite{quarteroni2017integrated}. 

For the space discretization, we use the Finite Element (FE) method with continuous FE on a tetrahedral mesh. We consider FE of order 2 ($\mathbb{P}_2$) for the electrophysiology to capture the traveling wave dynamics and FE of order 1 ($\mathbb{P}_1$) for both the activation and the mechanics~\cite{fedele2022comprehensive}.

For the time discretization, we use finite difference schemes~\cite{quarteroni2009numerical}. Specifically, cardiac electrophysiology is solved with Backward Differentiation Formula (BDF) of order 2, using an Implicit-Explicit (IMEX) scheme where the diffusion term is treated implicitly, the ionic and reaction terms explicitly. The ionic variables are advanced in time through an IMEX scheme~\cite{regazzoni2022cardiac,piersanti2022biv,fedele2022comprehensive}. We solve the active contraction problem with an IMEX BDF1 method, and the mechanical problem with a fully-implicit BDF1 scheme~\cite{piersanti2022biv}. Finally, an IMEX scheme of the first order is used for the circulation~\cite{fedele2022comprehensive}. Moreover, two different time steps are used: a finer one for the electrophysiology and a larger one for both the activation, the mechanics and the circulation~\cite{regazzoni2022cardiac}. Finally, we employ recently developed stabilization methods – related to the circulation and the fibers-stretch-rate feedback – that are crucial to obtain a stable solution in a four-chamber simulation scenario~\cite{regazzoni2021oscillation,regazzoni2022stabilization}. 
Concerning the linear systems arising from the discretization of the whole-heart EM problem we use: the conjugate gradient for the electrophysiology and the GMRES method for both the mechanics and the activation, {both empowered by} an {algebraic multigrid (AMG)} preconditioner.    Finally, we solve the non-linear saddle-point problem arising from the coupling between the mechanics and the circulation by means of a Newton algorithm using, at the algebraic level, the Schur complement reduction ~\cite{piersanti2022biv, regazzoni2022cardiac}.

For further details about the numerical methods we use in the whole-heart EM model, we refer to~\cite{regazzoni2022cardiac,piersanti2022biv,fedele2022comprehensive}.

\paragraph{Numerical methods for the {fluid domain displacement} problem}
After simulating the whole-heart EM and reaching a limit cycle in terms of pressure and volumes, we extract the solution {from} the last simulated heartbeat. We restrict {this} solution to the heart endocardium and the endothelium of outflow tracts, obtaining $N + 1$ solutions defined on the boundary of the fluid domain ($\bm d^{\partial \Omega}_i$, with $i = 0, \dots, N$). We project $\bm d^{\partial \Omega}_i$ onto the CFD mesh {with piecewise linear interpolation}, then solve the {fluid domain displacement} problem in \Cref{eq:lifting} to {obtain} the ALE displacement $\bm d$. 
We discretize the lifting problem in \Cref{eq:lifting} using FEs of order 1 ($\mathbb P_1$) and the linear system arising from its discretization is preconditioned with an {AMG} preconditioner.  We solve the resulting linear system with the conjugate gradient method. The smoothing spline approximation is computed independently for each mesh node, and the approximant is constructed following the optimization procedure described in~\cite{rodriguez2001smoothing}. To compute the ALE velocity, we use BDF1 to discretize in time \Cref{eq:velocity_ALE}.

\paragraph{Numerical methods for the CFD cardiac model}
We discretize \Cref{eq:ns_ale_riis} in space {using} FEs of order 1 for both velocity and pressure ($\mathbb P_1-\mathbb P_1$). We employ the Variational Multiscale - Large Eddy Simulation (VMS-LES) method to obtain a stable formulation of the NS-ALE-RIIS equations discretized via equal order FE spaces. This also allows us to control instabilities arising from the advection-dominated regime, and to model transition-to-turbulence in the LES framework~\cite{forti2015semi, bazilevs2007variational, zingaro2021hemodynamics}. The VMS-LES formulation accounts for the ALE framework and the RIIS modeling used for valves. For the complete formulation, we refer to~\cite{zingaro2022geometric}.

For the time discretization, we consider a uniform partition of the temporal domain in $N_t$ subintervals $(t_n, t_{n+1}]$ of uniform size $\Delta t$, with $n = 0, \dots, N_t-1$. We denote from here quantities approximated at time $t_n$ with the subscript $n$, e.g. $\bm u_n \approx \bm u(t_n)$. We advance the problem in time by means of BDF1. To reduce the computational burden of the numerical simulations, we use a semi-implicit treatment of the nonlinearities, as done in~\cite{forti2015semi}. The overall numerical scheme for the fluid dynamics problem is detailed in~\cite{zingaro2022geometric}.

Since Neumann boundary conditions may give rise to instability phenomena in case of inflow, we set backflow stabilization on all the Neumann boundaries in the inertial form presented in~\cite{bertoglio2014tangential}. 

The linear system arising from the discretization of \Cref{eq:ns_ale_riis} is preconditioned with the aSIMPLE preconditioner~\cite{Deparis_2014}, and each of its blocks is preconditioned with an AMG preconditioner. The linear system is then solved at each time step {by} the GMRES method.

\paragraph{Numerical method for the 0D circulation model}
We solve the system of ODEs of the circulation problem with an IMEX method of the first order. The time-step size $\Delta t$ employed for its numerical discretization is the same used for the BDF advancing scheme in the 3D problem.

\paragraph{Numerical scheme for the coupled CFD problem}

After initialization, for each temporal step of the CFD problem, we update the position of valve leaflets and we compute the ALE velocity. At every time step, we solve the 3D and 0D subproblems independently. First, we solve the 0D open circulation problem using $Q^\mathcal{I}_n$ as input (i.e. the flowrates computed in the 3D model at previous timestep, namely  $Q^\mathrm{SYS}_{\mathrm{VEN}, n}, \, Q_{\mathrm{PV}, n}, \, Q^\mathrm{PUL}_{\mathrm{VEN}, n}, \, Q_{\mathrm{AV}, n}$). Then, from the solution of the circulation, we compute the pressures $p^\mathcal{I}_{n+1}$ at the interfaces and solve the fluid dynamics problem providing those pressures as Neumann boundary conditions at inlet and outlet sections. Finally, we compute the interface data from the 3D to the 0D model, i.e. $Q^\mathcal{I}_{n+1}$. This approach treats the coupling between the 3D and 0D subproblems in a segregated and explicit way.


\section{Numerical results}
\label{sec:numerical-results}

In this section, we present the numerical results using the whole-heart fluid dynamics model. In \Cref{sec:setup-em-cfd}, we introduce the computational setting of the whole-heart EM and CFD simulations. \Cref{sec:calibration-rdq20mf} is devoted to the \rev{manual} calibration of the RDQ20 \rev{active contraction model} to produce physiological flowrates in the EM simulation. The physiological results of the overall computational model are presented in \Cref{sec:results-physiology}. Finally, we apply the multiphysics computational model to the pathological case of LBBB in \Cref{sec:application}. 

\subsection{Computational setup}
\label{sec:setup-em-cfd}

\begin{figure}[t]
	\centering
	\includegraphics[width=\textwidth]{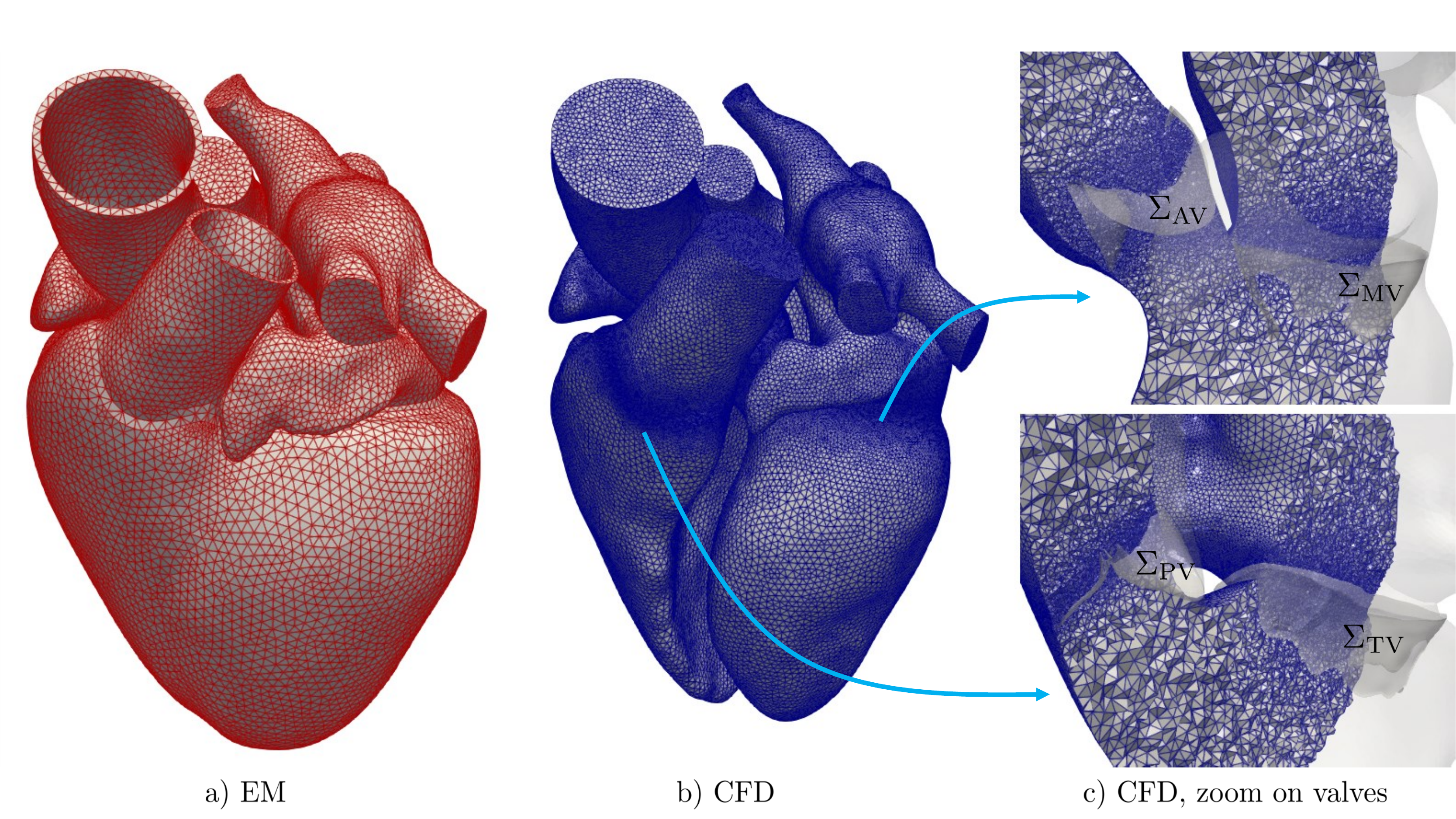}
	\caption{Tetrahedral meshes used in the computational model: a) mesh for the EM simulation, b) mesh for the CFD simulation, c) mesh refinement on the valves region of the CFD mesh.}
	\label{fig:mesh-em-cfd-zoom}
\end{figure}

We consider a realistic whole-heart geometry provided by the Zygote solid 3D heart model~\cite{zygote}, an anatomically CAD model representing an average healthy human heart reconstructed from high-resolution computer tomography scan data. We generate whole-heart tetrahedral meshes for the EM and CFD problems that we report in \Cref{fig:mesh-em-cfd-zoom}. Meshes are generated with \texttt{vmtk}~\cite{vmtk} using the methods and tools discussed in~\cite{fedele2021polygonal, fedele2022comprehensive, zingaro2022geometric}. Details on the meshes for the EM and CFD simulations are provided in \Cref{tab:mesh-em}.
The valve leaflets are thin structures that we characterize, in the context of the RIIS method, by small values of $\varepsilon_{\mathrm k}$. To correctly capture the immersed surfaces, we refine the CFD mesh close to the valve regions, as shown in \Cref{fig:mesh-em-cfd-zoom}c. Specifically, following~\cite{fedele2017patient}, we choose $h_\text{k}$ such that $\varepsilon_\mathrm{k} \geq 1.5 h_\text{k}$, where $h_\text{k}$ is the minimum mesh size of the fluid mesh in the valve region. 

\begin{table}[t!]
	\centering
	\begin{tabular}{ccccccccc}
		\toprule
		\textbf{Simulation} & \multicolumn{3}{c}{\textbf{Mesh size} {[\si{\milli\metre}]}} & {\textbf{Cells}} &  {\textbf{Points}} & {\textbf{Physics}} & {\textbf{DOFs}}  & {$\Delta t$} [\si{\second}]
		\\
		& \footnotesize{min} & \footnotesize{avg} & \footnotesize{max}  & & &  & \\
		\midrule
		\multirow{3}{*}{\textbf{EM}} & \multirow{3}{*}{\num{0.860}} & \multirow{3}{*}{\num{2.97}} & \multirow{3}{*}{\num{5.34}} & \multirow{3}{*}{\num{180472}}  &  \multirow{3}{*}{\num{46915}}  & Electrophysiology & \num{310505}& \num{5e-5}
		\\
		& &  &  &  &  & Mechanics  & \num{140745}& $10^{-3}$
		\\
		& &  &  &  &  & Circulation  & - &$10^{-3}$
		\\
		\midrule
		\multirow{2}{*}{\textbf{CFD}} & \multirow{2}{*}{\num{0.210}} & \multirow{2}{*}{\num{1.04}} & \multirow{2}{*}{\num{3.82}} & \multirow{2}{*}{\num{3892584}}  &   \multirow{2}{*}{\num{652204}} & Fluid dynamics & \num{2608816} & $10^{-4}$
		\\
		&  &  &  &  &  &  Circulation & - &$10^{-4}$
		\\
		\bottomrule
	\end{tabular}
	\caption{Setup of whole-heart EM and CFD simulations (mesh details and time step sizes).}
	\label{tab:mesh-em}
\end{table}

We carry out EM and CFD simulations in \texttt{life$^\texttt{x}$}~\cite{AFRICA2022101252, africa2022lifex, africa2024lifexcfd, africa2023lifexep}\footnote{\url{https://lifex.gitlab.io/}}, a high-performance \texttt{C++} FE library developed within the iHEART
project\footnote{iHEART - An Integrated Heart model for the simulation of the cardiac function, European Research Council (ERC) grant agreement No 740132, P.I. A. Quarteroni, 2017-2023}, mainly focused on cardiac simulations and based on the \texttt{deal.II} finite element core~\cite{arndt2021dealii,arndt2020dealii,dealii}.  

Numerical simulations are run in parallel on the GALILEO100 supercomputer\footnote{528 computing nodes each 2 x CPU Intel CascadeLake 8260, with 24 cores each, 2.4 GHz, 384GB RAM. See \url{https://wiki.u-gov.it/confluence/display/SCAIUS/UG3.3\%3A+GALILEO100+UserGuide} for technical specifications.} at the CINECA supercomputing center, using 240 and 480 cores for the EM and CFD simulations, respectively. The computational time to carry out a single heartbeat is about 1 hour and 20 minutes for the EM simulation and 56 hours for the CFD simulation. 

For the parameters of the EM model, we use the same values as in~\cite{fedele2022comprehensive} with some minor differences reported in \ref{appendix:em-cfd}. {As we better discuss in \Cref{sec:calibration-rdq20mf}, a huge difference in terms of setup of the EM simulation between~\cite{fedele2022comprehensive} and the present work consists in the \rev{manual} calibration of the \rev{active contraction model} to compute physiological blood flowrates}.  We simulate 20 heartbeats of the whole-heart EM, and we report numerical results related to the last heartbeat, after verifying that the solution is sufficiently close to a periodic limit cycle (in terms of pressure and volume transients). We consider an heartbeat period of $T_\mathrm{HB}=\num{0.8} \, \si{\second}$. 
We pick the last simulated EM heartbeat as input displacement for the CFD simulation. We set as initial condition for the velocity $\bm u_0 = \bm 0$. The initial state of the circulation model (for the coupling with the fluid dynamics) is taken equal to the values reached at the beginning of the last heartbeat in the EM simulation. Moreover, as detailed in \ref{appendix:em-cfd}, all the values of the parameters involved in the circulation model are {the same for} the EM and the fluid dynamics simulations. The physical parameters for blood are density {$\rho = \SI{1.06e3}{\kilo\gram\per\metre\cubed} $} and dynamic viscosity $\mu=\SI{3.5e-3}{\kilo\gram\per\metre\per\second}$.
We simulate two heart cycles, and we report the solution on the second cycle to remove the consequences of an unphysical null initial condition. For the numerical results visualization (for both EM and CFD), we shift the time domain in $(0, T_\mathrm{HB})$.

\begin{figure}[t]
	\centering
	\includegraphics[width=\textwidth]{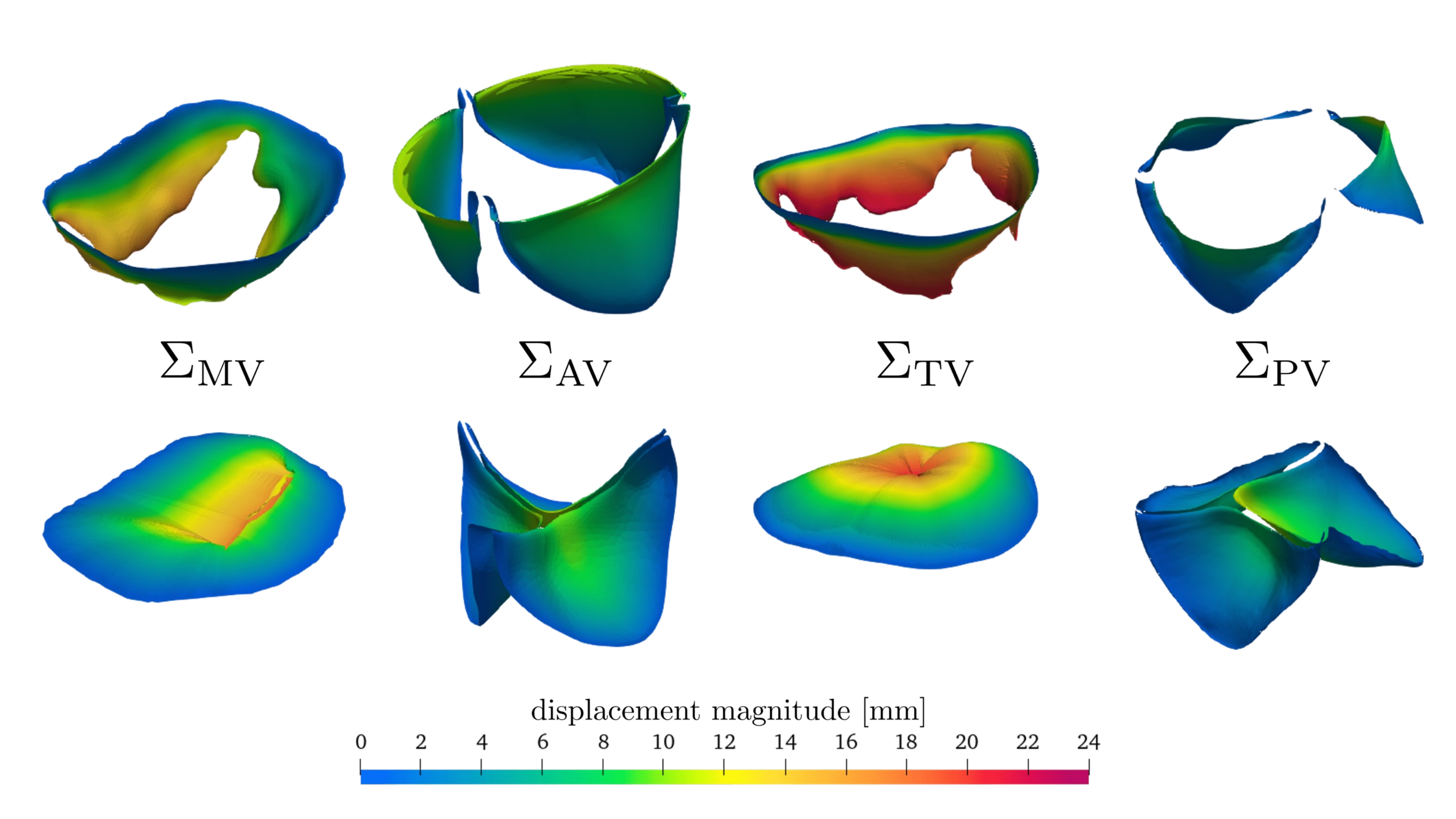}
	\caption{Cardiac valves in their open and closed configurations coloured according to displacement magnitude. Valves geometry are provided by Zygote~\cite{zygote}, and we define displacement field aimed at closing and opening the leaflets.}
	\label{fig:valves}
\end{figure}

The Zygote cardiac valves~\cite{zygote} are provided in their open configuration (TV, MV) and closed configuration (PV, AV). Thus, we define displacement fields aimed at closing and opening their leaflets, based on signed-distance functions and the solution of Laplace-Beltrami problems~\cite{fedele2021polygonal, zingaro2022geometric}. In \Cref{fig:valves}, we report the valves in their open and closed configurations, colored according to the leaflets' displacement magnitude.  We open and close the valves instantaneously (i.e. in one time step) at the times reported in \Cref{tab:valves-riis-setup}. These times are chosen by selecting the initial and final times of the isovolumetric phases. The values we choose for $\varepsilon_\mathrm{k}$,  reported in \Cref{tab:valves-riis-setup}, allow to have a physiological representation of the valve leaflet. Indeed, we choose $\varepsilon_\mathrm{k}$ by averaging the values of the leaflet thicknesses reported in~\cite{CRAWFORD20011419}. Furthermore, the valve resistances values $R_\mathrm{k}$ are reported in \Cref{tab:valves-riis-setup}. We found that the condition number of the linear system associated to the FE discretization of the fluid dynamics problem becomes larger as the ratio $R_\mathrm{k}/\varepsilon_\mathrm{k}$ increases. Thus, to keep contained the computational cost of the CFD simulation, we choose as $R_\mathrm{k}$ the minimum value that guarantees impervious valves.

\begin{table}[t]
	\centering
	\begin{tabular}{c c c c c c}
		\toprule
		& & \textbf{MV} & \textbf{AV} & \textbf{TV} & \textbf{PV} 
		\\  
		\midrule
		opening time & [\si{\second}] & \num{0.710} & \num{0.262} & \num{0.700}  & \num{0.279}
		\\
		closing time & [\si{\second}] & \num{0.208} & \num{0.666} & \num{0.194} & \num{0.677}
		\\
		$R_{\mathrm k}$ & [\si{\kilo\gram\per{\metre\per\second}}] & $10^{4}$ & $10^{4}$ &$10^{4}$ & $10^{4}$
		\\
		$\varepsilon_\mathrm k$ & [\si{\milli\metre}] & \num{0.68} & \num{0.67} & \num{0.77} & \num{0.52} 
		\\
		\bottomrule
	\end{tabular}
	\caption{Setup of the RIIS method to model cardiac valves. 
	}
	\label{tab:valves-riis-setup}
\end{table}


\subsection{Calibration of the RDQ20 \rev{active contraction model} to achieve physiological flows}
\label{sec:calibration-rdq20mf}

Among the different components at the basis of cardiac EM simulations, the model describing the force generation at the microscale plays a pivotal role. Indeed, not only the amount of force the muscle develops depends on it, but also its temporal distribution over the heartbeat, i.e. the kinetics of contraction and relaxation. \rev{As demonstrated in \cite{fedele2022comprehensive}, a biophysically detailed force generation model can effectively reproduce flowrates through the semilunar valves, which are very often severely overestimated in models involving less detailed active contraction models \cite{gerach2021electro, feng2024whole}}. Since the electromechanical displacement drives the fluid dynamics model, \rev{the CFD flow rates are the same as those of the electromechanical simulation, and are therefore heavily affected by the choice of the force generation model. Thus,} special care must be devoted to the choice and calibration of the \rev{active contraction model}.

For the above reasons, we chose to use the RDQ20 model~\cite{regazzoni2020biophysically}, that is an active force model with high biophysical fidelity and that is able to reproduce the main {features} of the experimentally observed behaviors.
The RDQ20 model is based on a detailed description of the calcium-driven regulation of the thin filament, with explicit representation of end-to-end cooperative interactions, and a description of the attachment-detachment process of crossbridges, at the basis of the force-velocity relationship.
Thereby, the model is able to reproduce the main mechanisms of contractility regulation, mediated by calcium, fiber strain and fiber strain-rate.
In particular, the fiber strain-rate feedback, which is responsible for the well-known force-velocity relationship, plays a central role in the regulation of hemodynamic flows, as demonstrated in~\cite{fedele2022comprehensive} and confirmed in the present study.

\newcommand{\RDQcmame}{A}
\newcommand{\RDQslowCa}{B}
\newcommand{\RDQslowCb}{C}
\newcommand{\RDQslowCc}{D}
\newcommand{\RDQsrf}{E}

On this basis, we refine the calibration of the RDQ20 model, with a particular care on  fluxes through semilunar valves obtained by means of the 0D model in the EM simulation. We employ as a starting point the calibration used in~\cite{fedele2022comprehensive}, suitable for the coupling with the TTP06 ionic model~\cite{TTP06} (see \Cref{tab:RDQ20-setup}, setting {\RDQcmame}). \rev{In this paper, we focus on the calibration of the ventricles only, whereas for the atria, we employed the same parameters of the RDQ20 model used in \cite{fedele2022comprehensive}.}
In \Cref{tab:RDQ20-results}, column {\RDQcmame}, we report a list of biomarkers obtained by using the calibration {\RDQcmame} in the EM simulation.  \rev{To better appreciate the effect of each calibration on different biomarkers, we provide a graphical representation in \Cref{fig:calibration-rdq20-graphic}. We normalize each biomarker in the interval $[-1, 1]$ using reference physiological values from the literature and reported in \Cref{tab:RDQ20-results}. }
Although the biomarkers characterizing the overall cardiac function (i.e. end-systolic and end-diastolic volume, stroke volume and ejection fraction) are within reference ranges, the maximum blood flux across valves is significantly above the physiological range.
In other terms, even if the total ejected blood is physiological, the instantaneous flow peak is too large.
This would clearly have a strong negative impact on the results of fluid dynamics simulations. For instance, an excessively high velocity through the valve may result in high pressure gradients, and an overall incorrect stress distribution over valve leaflets and cardiac walls \cite{van2017aortic, kriz2010renal}.

\newcommand{\modifiedparam}[1]{(*)\textbf{#1}}

\begin{table}
	\centering
	\begin{tabular}{c S c c c c c}
		\toprule
		\multicolumn{2}{l} {\textbf{Parameter}} & 
		\RDQcmame & \RDQslowCa & \RDQslowCb & \RDQslowCc & \RDQsrf
		\\
		\midrule
		\multicolumn{7}{l} {\textbf{Regulatory units dynamics}}  \\
		$Q$                               & [\si{-}]                            & 2       & 2       & 2       & 2       & 2 \\
		$\overline{k}_\mathrm{d}$         & [\si{\micro\molar}]                 & 0.36    & 0.36    & 0.36    & 0.36    & 0.36 \\
		$\alpha_{k_\mathrm{d}}$           & [\si{\micro\molar\per\micro\meter}] & -0.2083 & -0.2083 & -0.2083 & -0.2083 & -0.2083 \\
		$\mu$                             & [\si{-}]                            & 10      & 10      & 10      &  10     & 10 \\
		$\gamma$                          & [\si{-}]                            & 30      & 30      & 30      &  30     & 30 \\
		$K_\mathrm{off}$                  & [\si{\per\second}]                  & 8       & \modifiedparam{4} & \modifiedparam{4} & \modifiedparam{4}  & \modifiedparam{4} \\
		$K_\mathrm{basic}$                & [\si{\per\second}]                  & 4       & \modifiedparam{2} & \modifiedparam{2} & \modifiedparam{2}  & \modifiedparam{2} \\
		\midrule[0.1pt]
		\multicolumn{7}{l} {\textbf{Crossbridge dynamics (prescribed)}}  \\
		$v_0$                             &  [\si{\per\second}]                 & 2 & 2 & 2 & 2 & \modifiedparam{0.5} \\
		$v_\mathrm{max}$                  &  [\si{\per\second}]                 & 8 & 8 & 8 & 8 & \modifiedparam{2}   \\
		$\tilde{k}_2$                     &  [\si{-}]                           & 66 & 66 & 66 & 66 & 66 \\
		$\mu_0^{\text{iso}}$              &  [\si{-}]                           & 0.22 & 0.22 & 0.22 & 0.22 & 0.22 \\
		\midrule[0.1pt]
		\multicolumn{7}{l} {\textbf{Crossbridge dynamics (automatically calibrated)}}  \\
		$r_0$                             &  [\si{\per\second}]                 & 134.31 & 134.31 & 134.31 & 134.31 & 33.24 \\
		$\alpha$                          &  [\si{-}]                           & 25.184 & 25.184 & 25.184 & 25.184 & 24.93 \\
		$\mu^0_{f_\mathcal P} $           &  [\si{\per\second}]                 & 32.225 & 32.225 & 32.225 & 32.225 & 7.98  \\
		$\mu^1_{f_\mathcal P} $           &  [\si{\per\second}]                 & 0.768  & 0.768  & 0.768  & 0.768  & 0.192 \\
		\midrule[0.1pt]
		\multicolumn{7}{l} {\textbf{Micro-macro upscaling}}  \\
		$a_\mathrm{XB}$                   &  [\si{\mega\pascal}]                & 1500.0 & \modifiedparam{1550.0} & \modifiedparam{2925.0} & \modifiedparam{5214.5} &  \modifiedparam{1550.0} \\
		\bottomrule
	\end{tabular}
	\caption{Different calibrations of the RDQ20 \rev{active contraction model}: 
		\RDQcmame) calibration from~\cite{fedele2022comprehensive};
		\RDQslowCa) \RDQslowCb) \RDQslowCc) Reduced $ K_\mathrm{off}, \, K_\mathrm{basic}$ and different levels of $a_\mathrm{XB}$;
		\RDQsrf) Reduced $ v_0, \, v_\mathrm{max}$.
		Parameters that are modified with respect to {\RDQcmame} are highlighted by the symbol \textbf{(*)}.
		We remark that the parameters $r_0$, $\alpha$, $\mu^0_{f_\mathcal P} $ and  $\mu^1_{f_\mathcal P} $ are automatically calibrated from the four quantities $v_0$, $v_\mathrm{max}$, $\tilde{k}_2$ and $\mu_0^{\text{iso}}$ (for furhter details, see~\cite{regazzoni2020biophysically}); therefore, the symbol \textbf{(*)} is not reported for these four parameters.
		For the description of each parameter, we refer to the original paper of the RDQ20 model~\cite{regazzoni2020biophysically}. }
	\label{tab:RDQ20-setup}
\end{table}

\newcommand{\bioOver}[1]{($\uparrow$)\textbf{#1}}
\newcommand{\bioBelow}[1]{($\downarrow$)\textbf{#1}}

\begin{table}
	\centering
	\begin{tabular}{c S c c r r r r r }
		\toprule
		\multicolumn{2}{c} {\textbf{Biomarker}} & 
		\multicolumn{2}{c} {\textbf{Physiological values}} &
		\RDQcmame & \RDQslowCa & \RDQslowCb & \RDQslowCc & \RDQsrf
		\\
		\midrule
		$\mathrm{ESV}_\mathrm{LV}$            & [\si{\milli\litre}]            &  35  to  80   &~\cite{maceira2006normalized}            & 53.8            &  66.7               & 53.5             & 45.9             & 66.4  \\ 
		$\mathrm{ESV}_\mathrm{RV}$            & [\si{\milli\litre}]            &  69 $\pm$ 22  &~\cite{hudsmith2005normal}               & 58.0            &  79.8               & 65.3             & 56.0             & 72.3  \\ 
		$\mathrm{EDV}_\mathrm{LV}$            & [\si{\milli\litre}]            & 126 to 208    &~\cite{maceira2006normalized}            & 150           &  128             & \bioBelow{108} & \bioBelow{87.9}  & 151 \\ 
		$\mathrm{EDV}_\mathrm{RV}$            & [\si{\milli\litre}]            & 144 $\pm$ 23  &~\cite{maceira2006reference}             & 153           &  152              & 134            & 119            & 159 \\ 
		$\mathrm{SV}_\mathrm{LV}$             & [\si{\milli\litre}]            & 81 to 137     &~\cite{maceira2006normalized}            & 96.3            &  \bioBelow{61.5}    & \bioBelow{54.9}  & \bioBelow{42.0}  & 84.9  \\ 
		$\mathrm{SV}_\mathrm{RV}$             & [\si{\milli\litre}]            & 94 $\pm$ 15   &~\cite{maceira2006reference}             & 95.4            &  \bioBelow{72.3}    & \bioBelow{69.1}  & \bioBelow{63.2}  & 87.1  \\ 
		$\mathrm{EF}_\mathrm{LV}$             & [\%]                           &  49 to 73     &~\cite{clay2006normal}                   & 64.2            &  \bioBelow{48.0}    & 50.2             & \bioBelow{47.8}  & 56.1  \\ 
		$\mathrm{EF}_\mathrm{RV}$             & [\%]                           & 53 $\pm$ 6    &~\cite{burkhardt2019right}               & \bioOver{62.2}  &  47.5               & 51.4             & 53.0             & 54.6  \\ 
		$\mathrm{Q}_\mathrm{AV}^\mathrm{max}$ & [\si{\milli\litre\per\second}] & 427 $\pm$ 129 &~\cite{hammermeister1974rate}            & \bioOver{697} &  327              & 347            & 304            & 399 \\ 
		$\mathrm{Q}_\mathrm{PV}^\mathrm{max}$ & [\si{\milli\litre\per\second}] & 427 $\pm$ 129 &~\cite{hammermeister1974rate}$\dagger$   & \bioOver{756} &  397              & 427            & 443            & 478 \\ 
		$\mathrm{p}_\mathrm{LV}^\mathrm{max}$ & [\si{\mmHg}]                   & 119 $\pm$ 13  &~\cite{sugimoto2017echocardiographic}    & \bioOver{154} &  \bioBelow{99.5}    & \bioBelow{93.2}  & \bioBelow{80.7}  & 120 \\ 
		$\mathrm{p}_\mathrm{RV}^\mathrm{max}$ & [\si{\mmHg}]                   & 35 $\pm$ 11   &~\cite{bishop1997clinical}               & 37.2            &  34.9               & 38.2             & 41.6             & 32.2  \\ 
		\bottomrule
	\end{tabular}
	\caption{Effects of different calibrations of the RDQ20 \rev{active contraction model} on mechanics and hemodynamics biomarkers.
		Biomarkers highlighted by the symbols \bioOver{} and \bioBelow{} lie outside the reference ranges, denoting values too large or too small, respectively, compared with reference ranges.
		$\dagger$: for the maximum PV flowrate, in absence of clinical ranges from literature, we consider the same normal values of the AV flowrate.			
		Each column corresponds to a different calibration as detailed in \Cref{tab:RDQ20-setup}. }
	\label{tab:RDQ20-results}
\end{table}

\begin{figure}
	\centering
	\includegraphics[trim={3cm 1 3cm 1 },clip,width=\textwidth]{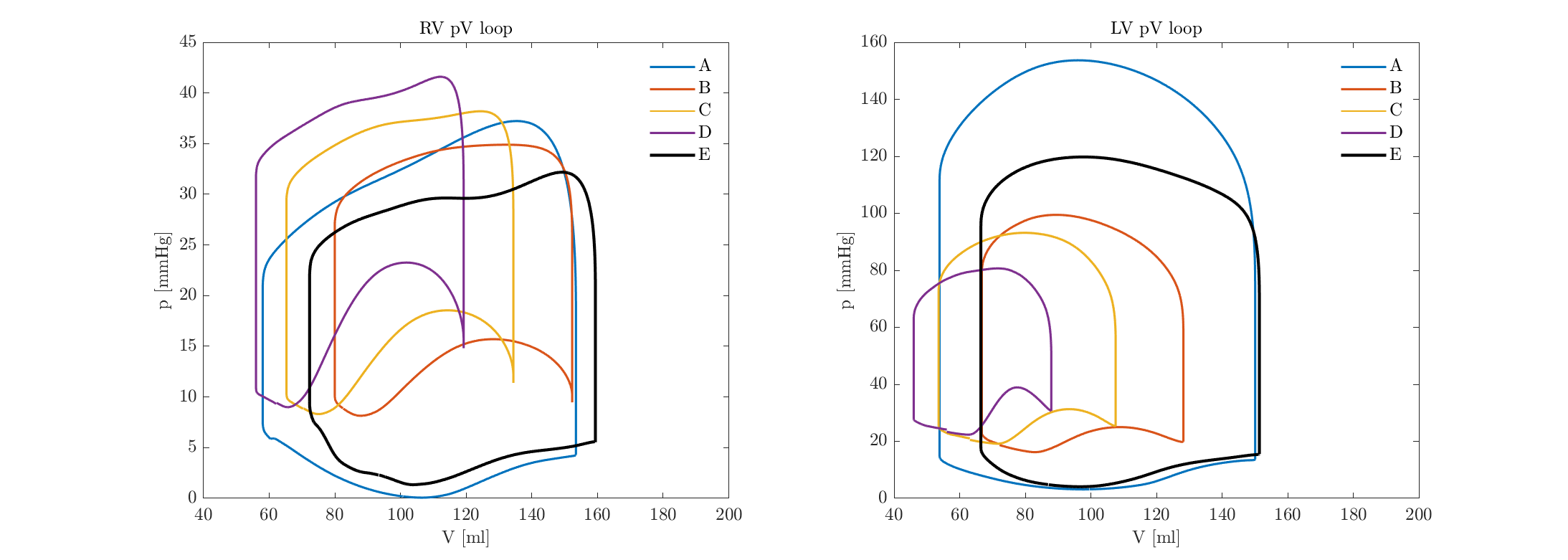}
	\caption{pV loops of RV and LV obtained with different calibrations of the RDQ20 \rev{active contraction model}. Each line corresponds to a different calibration as detailed in \Cref{tab:RDQ20-setup}.}
	\label{fig:pvloops-activation}
\end{figure}

To address this issue, we slow down the process of force generation, so that the tissue contractility is developed at a lower rate.
More precisely, we reduce the association-dissociation rates of troponin and tropomyosin of the RDQ20 model (i.e. $K_\mathrm{off}$ and $K_\mathrm{basic}$) by a factor 2. Moreover, to compensate for the lower peak force caused by a slower kinetics, we increase the crossbridge level contractility (i.e. $a_\mathrm{XB}$). 
We consider three different levels of contractility, as reported in \Cref{tab:RDQ20-setup}, respectively in columns {\RDQslowCa}, {\RDQslowCb} and {\RDQslowCc}.
However, as we show in \Cref{tab:RDQ20-results} and \Cref{fig:pvloops-activation}, on the one hand, we achieve the desired effect of reducing semilunar peak flows, thus bringing them {within} the expected ranges. On the other hand, we compute significantly reduced stroke volume and ejection fraction for both chambers, thus moving {out of} physiological ranges. Notice also that this issue is also not resolved by adjusting the contractility.
In fact, by raising $a_\mathrm{XB}$, not only the state of contractility in systole is changed, but also in diastole, leading to a reduced end-diastolic volume and, therefore, nullifying the effect of increased contractility, due to the Frank-Starling mechanism.
The three cases shown in \Cref{tab:RDQ20-results} ({\RDQslowCa}, {\RDQslowCb} and {\RDQslowCc}) are three illustrative cases out of the many that we tested, but without being able to reduce peak flows within the expected ranges while maintaining a physiological ejection fraction.
The tests evidenced a paradigmatic short-blanket problem, whereby just acting on kinetics and contractility it is not possible to lower the flows while maintaining a regular ejection fraction. 
Evidently, another element must be taken into account.

\begin{figure}
	\centering
	\begin{subfigure}{0.32\textwidth}
		\centering
		\includegraphics[trim={1 1 1 1 },clip,width=\textwidth]{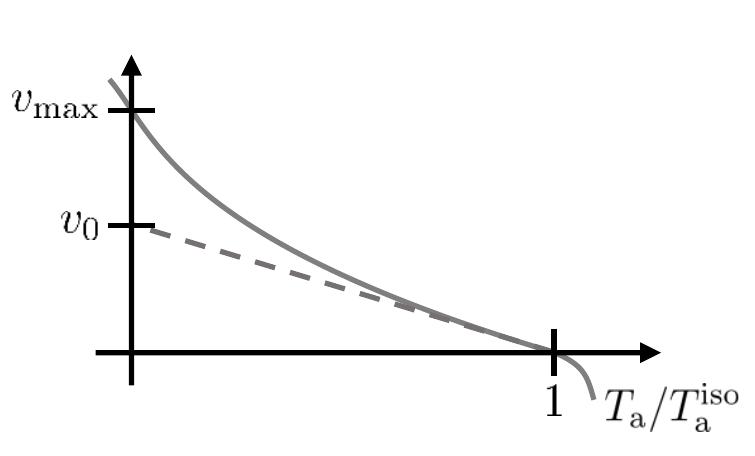}
		\caption{Force-velocity relation}
		\label{fig:calibration_F-vel}
	\end{subfigure}
	\begin{subfigure}{0.32\textwidth}
		\centering
		\includegraphics[trim={1 1 1 1 },clip,width=\textwidth]{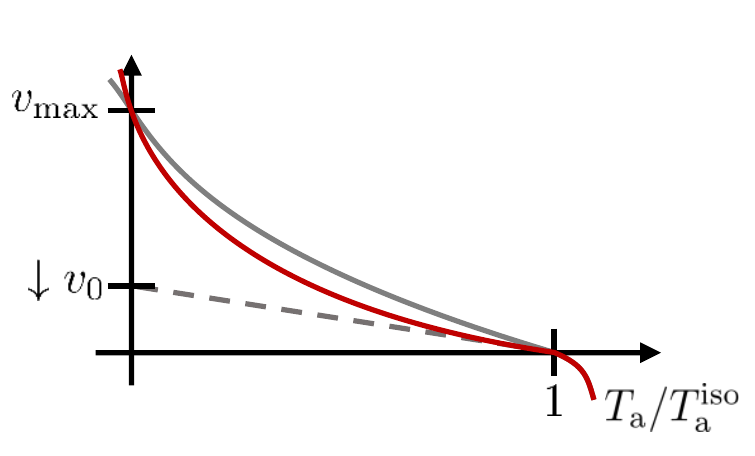}
		\caption{reduced $v_0$}
		\label{fig:calibration_F-vel-reducedv0}
	\end{subfigure}
	\begin{subfigure}{0.32\textwidth}
		\centering
		\includegraphics[trim={1 1 1 1 },clip,width=\textwidth]{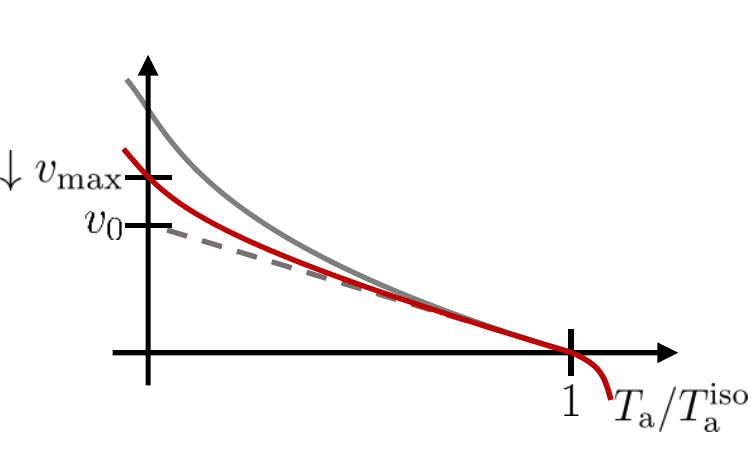}
		\caption{reduced $v_\mathrm{max}$}
		\label{fig:calibration_F-vel-reducedvmax}
	\end{subfigure}
	\caption{Representation of the well-known force-velocity relationship of muscle cells. The normalized active tension $T_{\mathrm{a}} / T_{\mathrm{a}}^{\mathrm{iso}}$, where $T_{\mathrm{a}}^{\mathrm{iso}}$ denotes the force in isometric conditions, is a decreasing function of the shortening velocity $v$.
		As shown in the figure, $v_\mathrm{max}$ is the shortening velocity for which the active tension reaches zero, while $v_0$ is the slope of the curve in correspondence of the isometric conditions (i.e. $v = 0$). (a) generic force-velocity relationship, (b) effect of reducing $v_0$ (in grey, the curve from (a)), (c) effect of reducing $v_\mathrm{max}$ (in grey, the curve from (a)).}
	\label{fig:calibration_F-vel-all}
\end{figure}

\begin{figure}[t]
    \centering
    \includegraphics[trim={0, 8cm, 0, 8cm},clip,width=\textwidth]{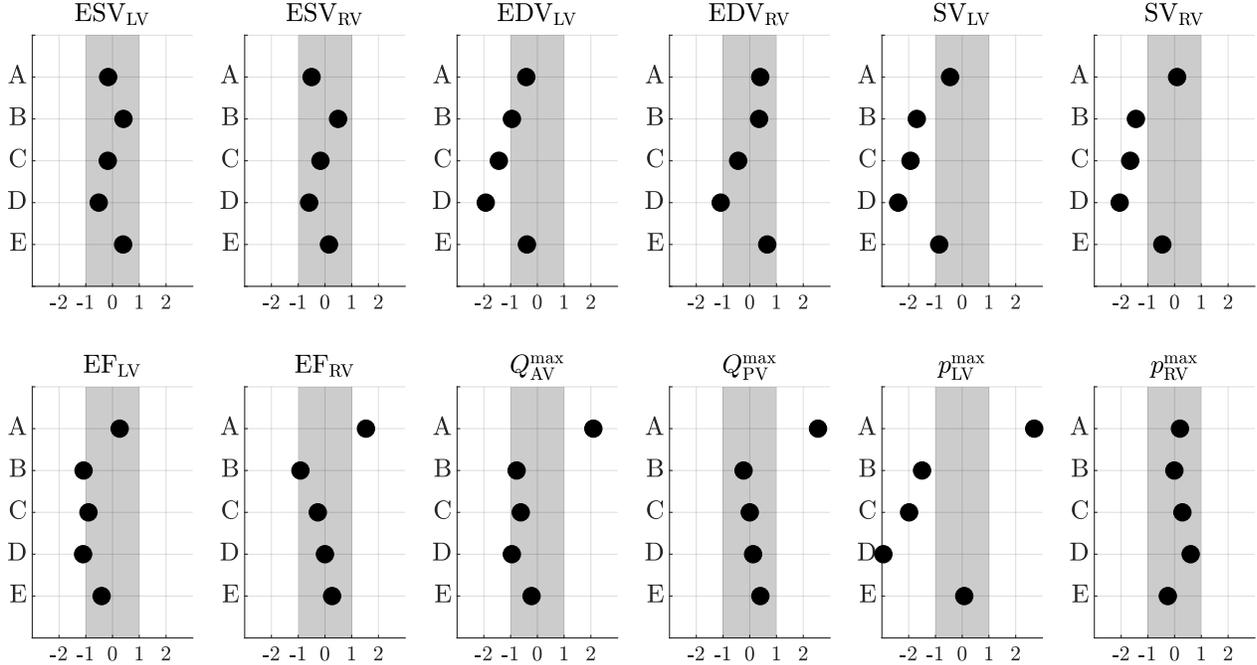}
    \caption{\rev{Effects of different calibrations of the RDQ20 active contraction model on mechanics and hemodynamics biomarkers. Each biomarker has been normalized according to the reference ranges provided in \Cref{tab:RDQ20-results} in the interval $[-1, 1]$. The gray rectangle denotes the reference range, whereas each black dots indicate the result of the EM simulation with a different calibration of the RDQ20 model, as detailed in \Cref{tab:RDQ20-setup}. }}
    \label{fig:calibration-rdq20-graphic}
\end{figure}

Based on the results of~\cite{fedele2022comprehensive}, which showed that, by neglecting the fibers-stretch-rate feedback, blood flows through the semilunar valves are significantly overestimated, we modified the calibration so as to, on the contrary, strengthen the effect of this feedback, but without changing either the isometric force or the kinetics.
Specifically, we acted in such a way as to steepen the force-velocity relationship, the microscopic mechanism underlying the fibers-stretch-rate feedback, by modifying the parameters governing cross-bridge dynamics.
For this purpose, we took advantage of the calibration technique illustrated in~\cite{regazzoni2020biophysically}, by which the RDQ20 model can be tuned to achieve a desired force-velocity relationship.
Two features of the force-velocity relationship can be selected, namely the maximum shortening velocity ($v_\mathrm{max}$), that is the velocity corresponding to vanishing active force, and the tangent to the curve under isometric conditions ($v_0$).
The geometric meaning of the two quantities is illustrated in \Cref{fig:calibration_F-vel-all}.

Hence, starting from the {\RDQslowCa} calibration, we modified the parameters to obtain a $v_\mathrm{max}$ and a $v_0$ equal to one-fourth of the original ones (see \Cref{tab:RDQ20-setup}, column {\RDQsrf}).
As evidenced in \Cref{tab:RDQ20-results}, with the setting {\RDQsrf} all biomarkers {fall} within physiological ranges (see also \Cref{fig:pvloops-activation}).
We believe that this result can be explained precisely by the mechanism of fibers-stretch-rate feedback, whereby regions of the tissue undergoing rapid shortening experience a decrease in developed force, thus promoting a more homogeneous shortening in space and without significant spikes in time, with a resulting viscous-like effect.
In conclusion, blood flow is redistributed more evenly over the duration of the ejection phase.

Based on the very good match with the reference values of the different biomarkers, in this work we {use} the calibration {\RDQsrf} as the baseline for EM simulations.


\subsection{Heart physiology and validation against clinical biomarkers}
\label{sec:results-physiology}

\begin{figure}
	\centering
	\begin{subfigure}{0.32\textwidth}
		\centering
		\includegraphics[trim={1 1 1 1 },clip,width=\textwidth]{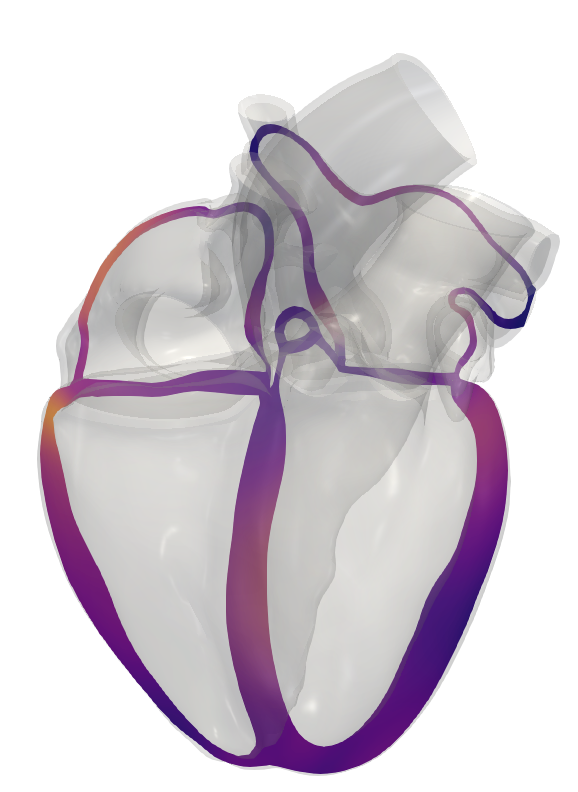}
		\caption{$t = 0.11$ s}
		\label{fig:em-displacement-0p11}
	\end{subfigure}
	\begin{subfigure}{0.32\textwidth}
		\centering
		\includegraphics[trim={1 1 1 1 },clip,width=\textwidth]{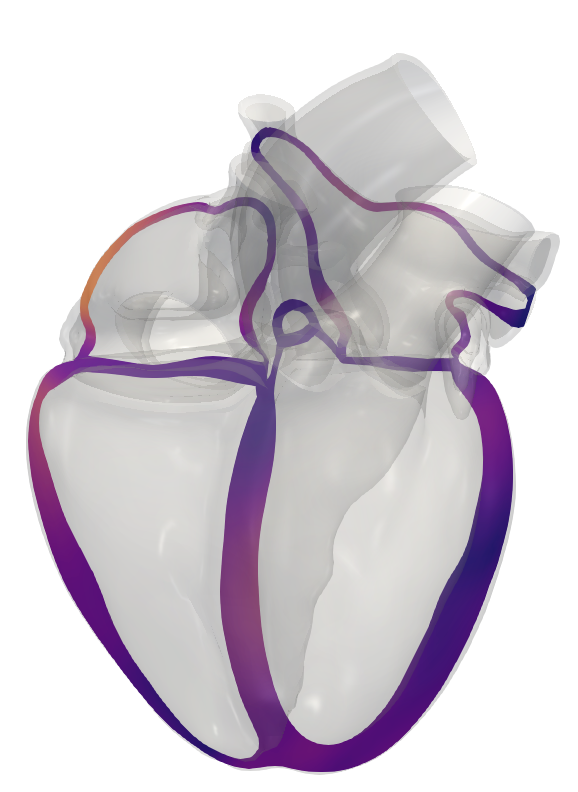}
		\caption{$t = 0.25$ s}
		\label{fig:em-displacement-0p25}
	\end{subfigure}
	\begin{subfigure}{0.32\textwidth}
		\centering
		\includegraphics[trim={1 1 1 1 },clip,width=\textwidth]{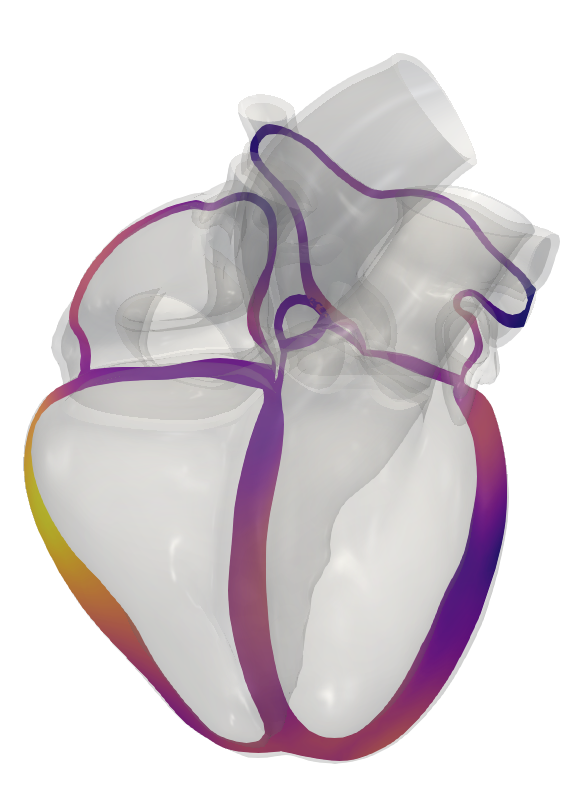}
		\caption{$t = 0.36$ s}
		\label{fig:em-displacement-0p36}
	\end{subfigure}
	\\
	\begin{subfigure}{0.32\textwidth}
		\centering
		\includegraphics[trim={1 1 1 1 },clip,width=\textwidth]{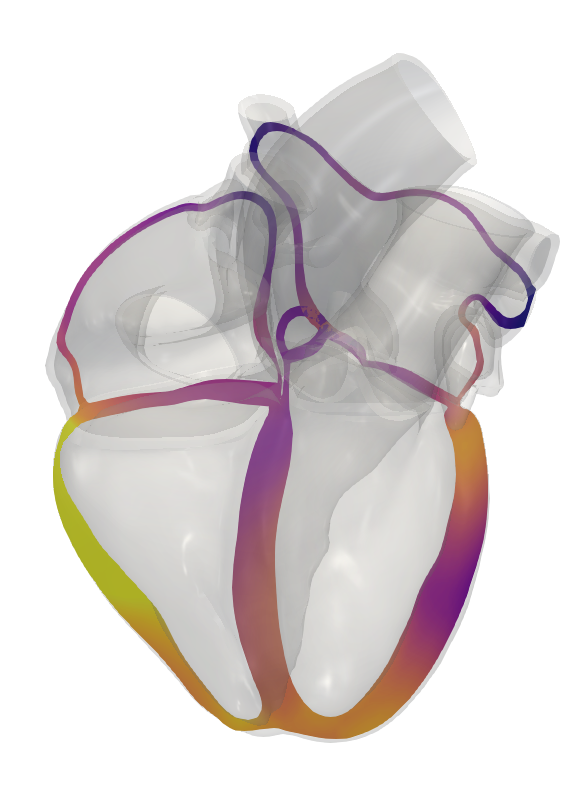}
		\caption{$t = 0.50$ s}
		\label{fig:em-displacement-0p50}
	\end{subfigure}
	\begin{subfigure}{0.32\textwidth}
		\centering
		\includegraphics[trim={1 1 1 1 },clip,width=\textwidth]{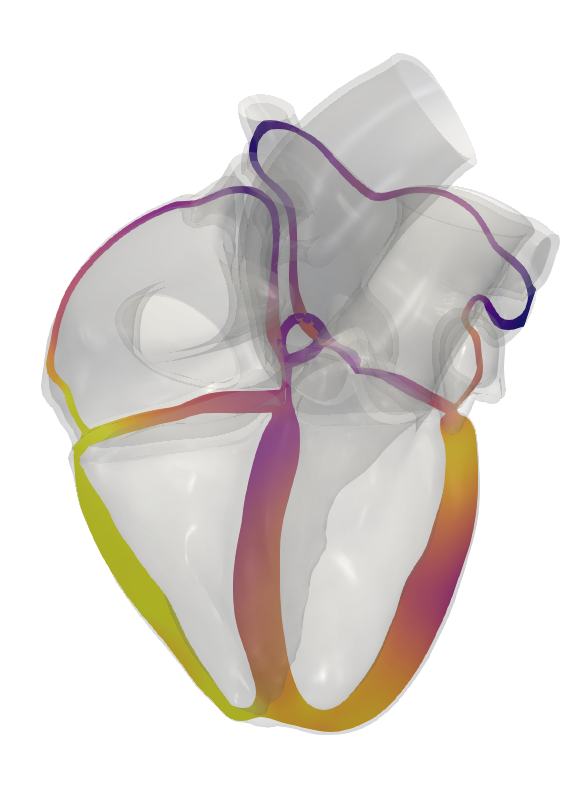}
		\caption{$t = 0.70$ s}
		\label{fig:em-displacement-0p70}
	\end{subfigure}
	\begin{subfigure}{0.32\textwidth}
		\centering
		\includegraphics[trim={1 1 1 1 },clip,width=\textwidth]{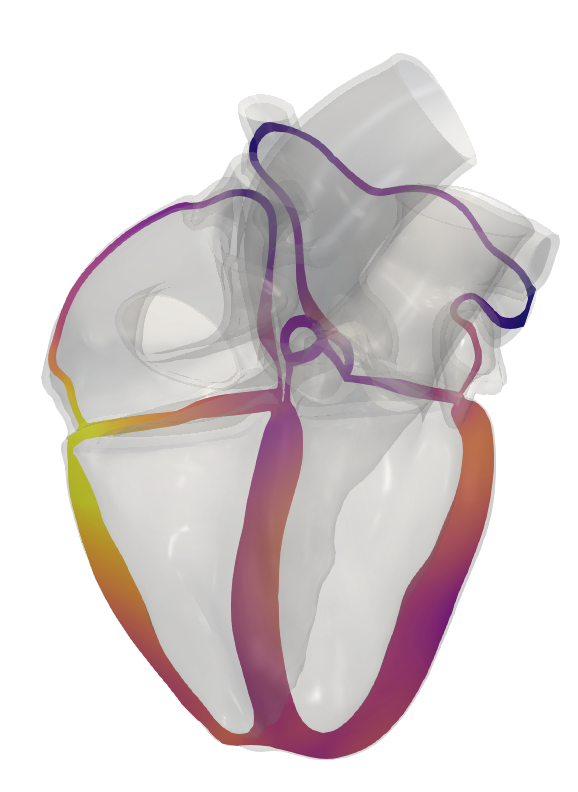}
		\caption{$t = 0.79$ s}
		\label{fig:em-displacement-0p79}
	\end{subfigure}
	\\
	\centering
	\vspace{0.5cm}
	$|\bm d^\mathrm{EM}| $ [mm]
	\\
	\centering
	\begin{subfigure}{\textwidth}
		\centering
		\includegraphics[trim={1 1 2cm 1 },clip,width=0.8\textwidth]{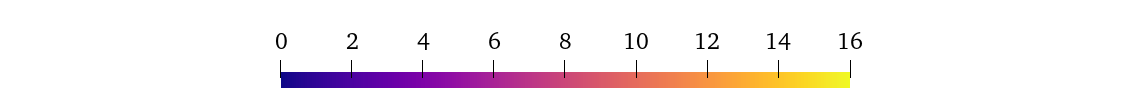}
	\end{subfigure}
	\caption{Whole heart deformed with EM displacement (with respect to the reference stress-free configuration) and colored according to its magnitude for a single, representative heartbeat: (a) active atrial contraction, (b) isovolumetric contraction, (c) ejection (peak), (d) ejection (mid-deceleration), (e) isovolumetric relaxation, (f) passive ventricular filling.}
	\label{fig:em-displacement}
\end{figure}

\Cref{fig:em-displacement} shows the whole-heart EM displacement for a single, representative heartbeat, where we highlight six different phases: isovolumetric contraction, ejection (peak and mid-deceleration), isovolumetric relaxation, ventricular passive filling, and atrial contraction. As pointed out in~\cite{fedele2022comprehensive}, the whole-heart EM model can correctly reproduce the cardiac physiological motion. Moreover, as shown in \Cref{sec:calibration-rdq20mf}, we found that, after a thorough calibration of the \rev{active contraction model}, our numerical results are consistent with normal values found in literature in terms of several volumetric biomarkers. 

\begin{figure}
	\centering
	\includegraphics[trim={3cm 1 3cm 1 },clip,width=\textwidth]{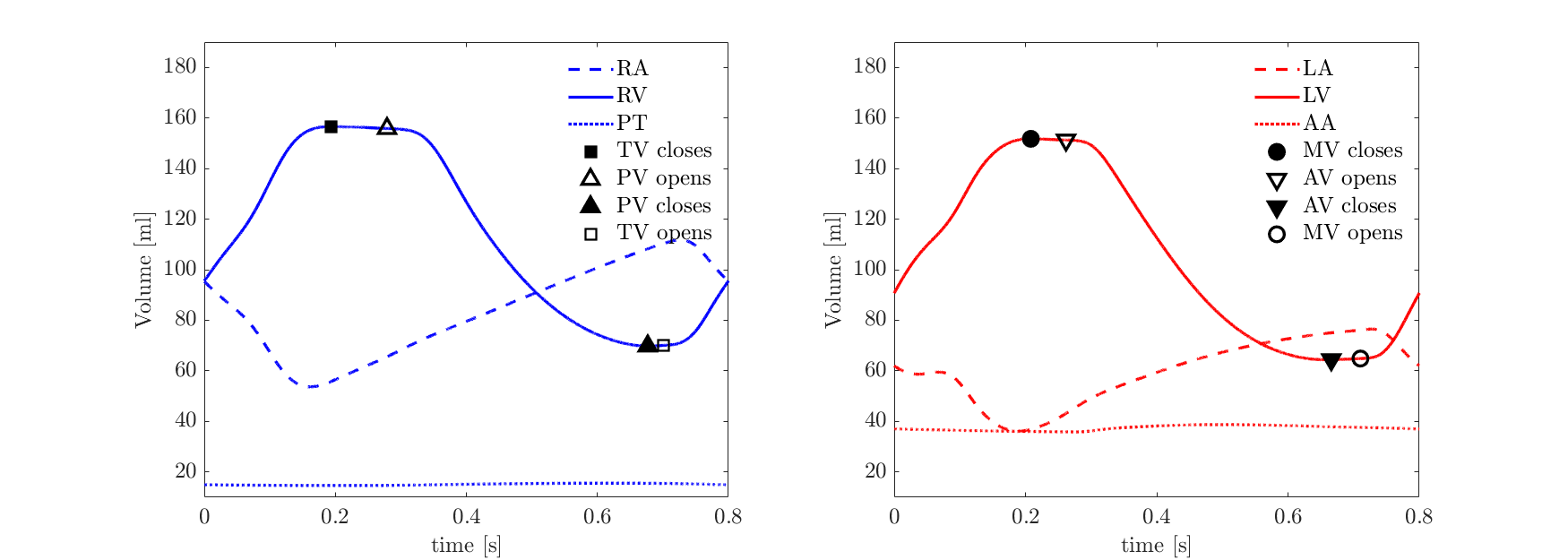}
	\caption{Volumes of RA, RV, PT (left) and LA, LV, AA (right) during a representative heartbeat. Valves open and close, instantaneously, at the initial and final times of isovolumetric phases. We report these times also in \Cref{tab:valves-riis-setup}.}
	\label{fig:volumes}
\end{figure}

In \Cref{fig:volumes}, we report the volumes of the heart chambers and large arteries versus time, and we highlight time instants in which valves open and close. In \Cref{fig:volume-rendering-velocity}, we report the volume rendering of the velocity magnitude obtained with our EM-driven CFD simulation. We start our {fluid dynamics} simulation at the end of ventricular diastole.

During the active atrial contraction (\Cref{fig:volume-rendering-0p10}), we observe the blood flowing from atria to ventricles, producing two high-speed jets in the MV and TV. This moment corresponds to the A-wave, as shown in \Cref{fig:velocities}. In order to assess whether {our} numerical simulation is correctly reproducing the physiological heart function, we compare the peak velocities through valves with physiological ranges available in literature and acquired in healthy subjects. We report this comparison in \Cref{tab:validation-cfd}. During diastole, we obtain lower velocities in the TV compared to MV, consistently with clinical measurements available in literature~\cite{thomas1998peak, choi2015new}. 

During the isovolumetric contraction (\Cref{fig:volume-rendering-0p24}), all valves are closed and the ventricular volumes remain constant. We measure lower velocity values compared to filling and ejection phases. Moreover, we found that the intraventricular pressure is not well defined and {is} prone to oscillations.  As a matter of fact, since we are using EM as unidirectional input of the CFD model, the dynamic balance between hemodynamics and tissue mechanics is neglected. Furthermore, since we are modeling the cardiac valves with the RIIS method, weakly imposing a kinematic condition only, the dynamic balance {is not fulfilled even} between blood and valves. Thus, our computational model cannot correctly capture the physiological pressure transient during this phase, instead producing nonphysical and large oscillations. In \Cref{fig:pressure}, we report the pressure transients in time and, for visualization purposes, we ignore the isovolumetric phases from the plot (grey boxes).

The ejection phase (\Cref{fig:volume-rendering-0p35} and \Cref{fig:volume-rendering-0p55}) is characterized by the opening of semilunar valves, the contraction of ventricles, and the blood flowing from the LV to the AO and from the RV to the PT. We measured peak flow rates equal to \num{398.74} \si{\milli\liter\per\second} and \num{478.16} \si{\milli\liter\per\second}, in the AV and PV, respectively, consistently with physiological values~\cite{hammermeister1974rate}. Moreover, as shown in \Cref{tab:validation-cfd}, we found that also maximum velocities between AV and PV and peak ventricular pressures during ejection are always in physiological ranges.  

During the isovolumetric relaxation, both the atrioventricular and semilunar valves are closed. The velocities measured are low compared to {those of} the other phases of the heartbeat. Furthermore, as for the isovolumetric contraction, the computational model cannot reproduce the typical pressure decrease occurring in this phase. On the contrary, large pressure oscillations arise. 

\begin{figure}
    \centering
    \includegraphics[trim={0, 9cm, 0, 9cm},clip,width=0.8\textwidth]{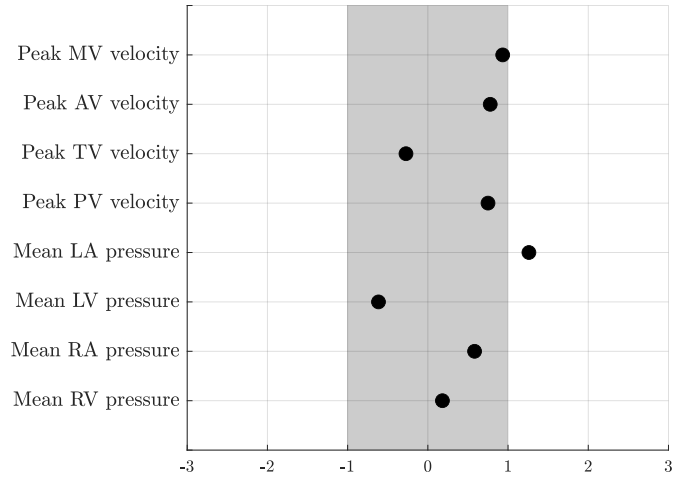}
    \caption{\rev{Fluid dynamics biomarkers obtained with the whole heart CFD simulations and normalized in the range $[-1, 1]$ with respect to the physiological ranges given in \Cref{tab:validation-cfd}. The grey rectangle denotes the reference range.}}
    \label{fig:cfd-normalized-graphic}
\end{figure}

\begin{figure}
	\centering
	\includegraphics[trim={3cm 1 3cm 1 },clip,width=\textwidth]{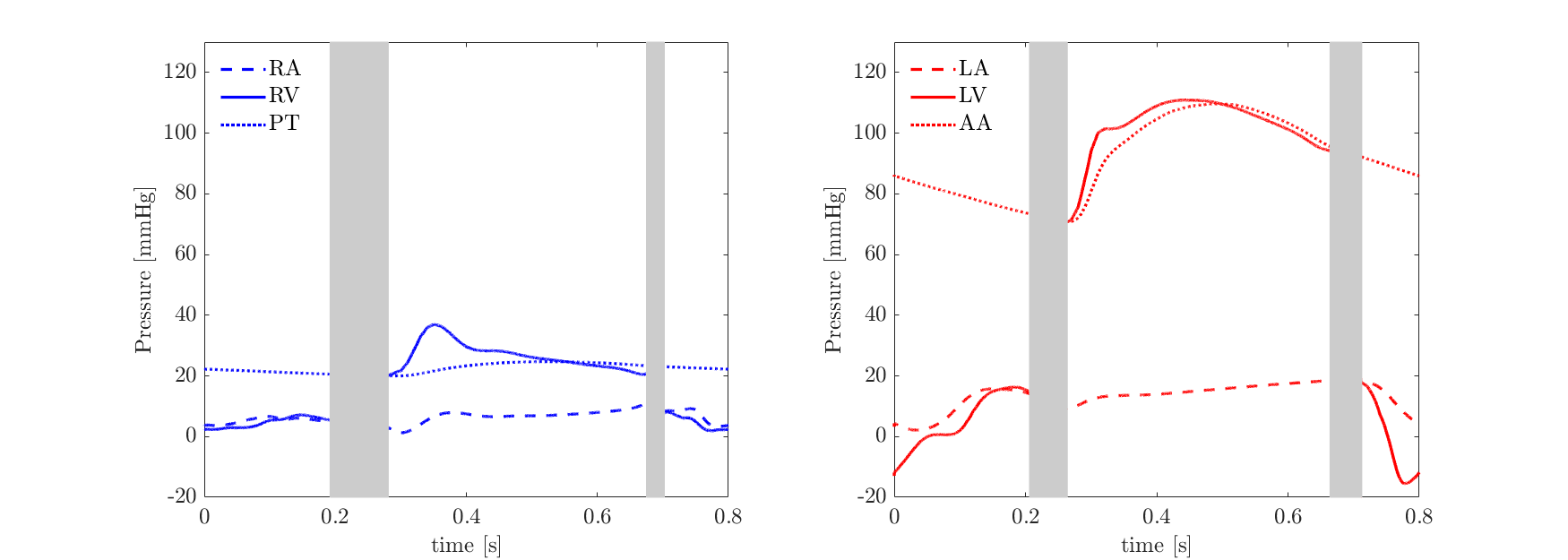}
	\caption{Pressure computed in control volumes located in RA, RV, PT (left) and LA, LV, AO (right) during a representative heartbeat. Isovolumetric phases are not considered, since the intraventricular pressure is not well defined when both valves are closed. }
	\label{fig:pressure}
\end{figure}

\begin{figure}
	\centering
	\includegraphics[trim={3cm 1 3cm 1 },clip,width=\textwidth]{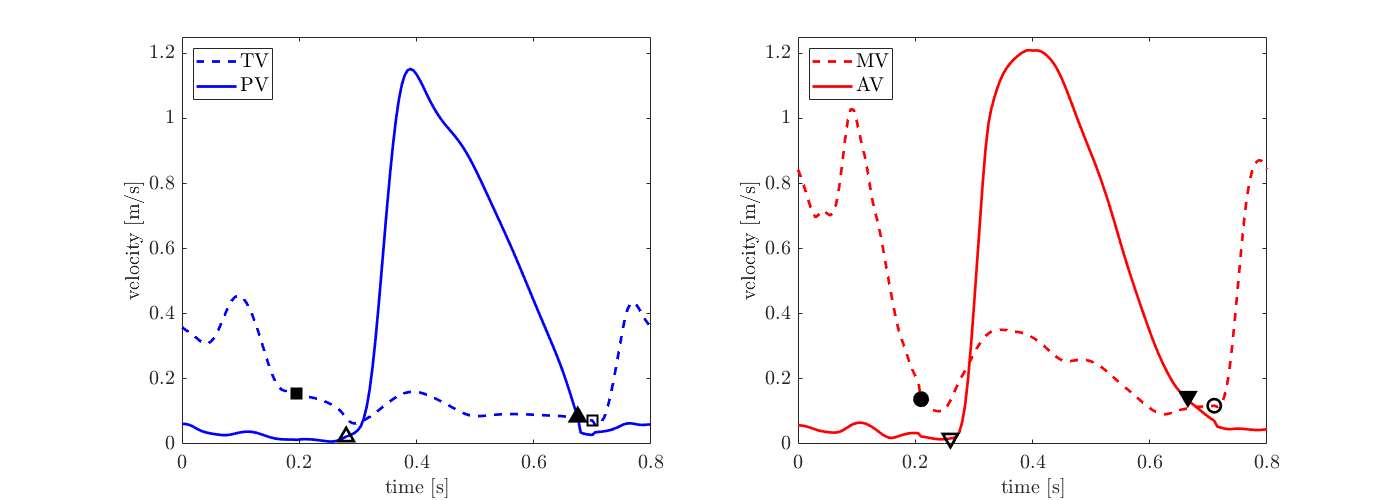}
	\caption{Velocity magnitudes computed in control volumes located between the valve leaflets: TV, PV (left) and MV, AV (right). }
	\label{fig:velocities}
\end{figure}

\begin{figure}[t]
	\centering
	\begin{subfigure}{0.32\textwidth}
		\centering
		\includegraphics[trim={1cm 1 1cm 1 },clip,width=\textwidth]{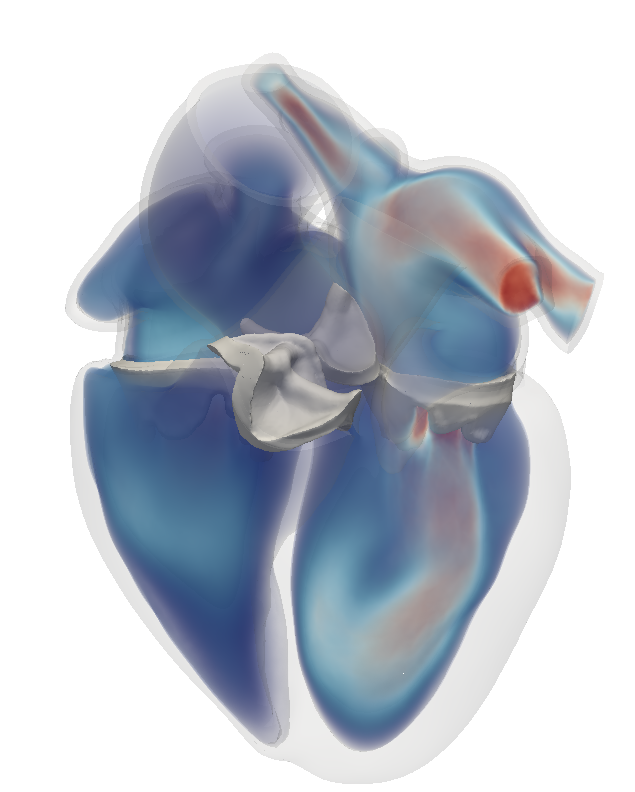}
		\caption{$t = 0.10$ s}
		\label{fig:volume-rendering-0p10}
	\end{subfigure}
	\begin{subfigure}{0.32\textwidth}
		\centering
		\includegraphics[trim={1cm 1 1cm 1 },clip,width=\textwidth]{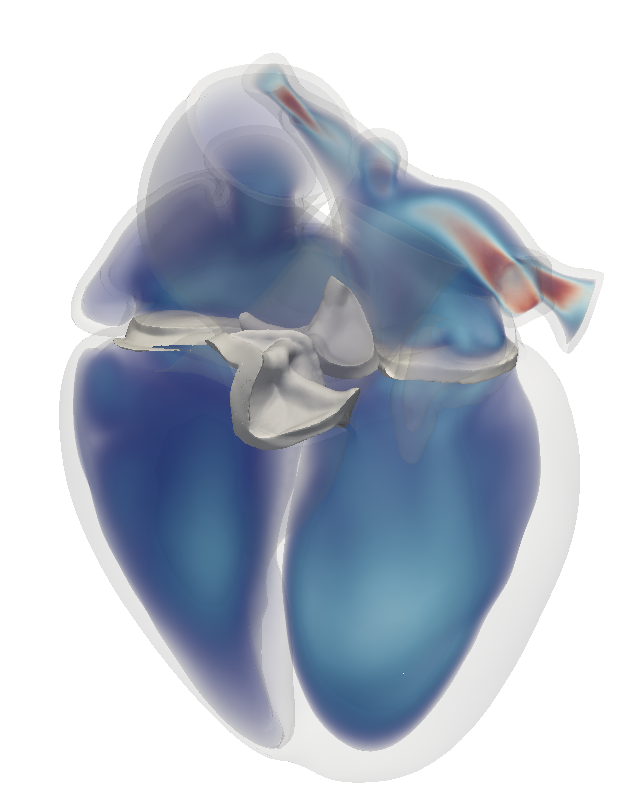}
		\caption{$t = 0.24$ s}
		\label{fig:volume-rendering-0p24}
	\end{subfigure}
	\begin{subfigure}{0.32\textwidth}
		\centering
		\includegraphics[trim={1cm 1 1cm 1 },clip,width=\textwidth]{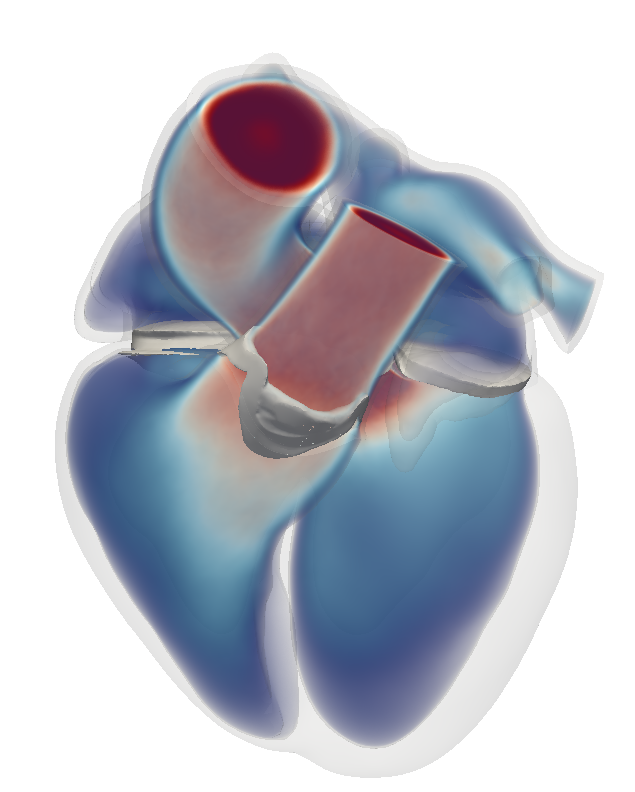}
		\caption{$t = 0.35$ s}
		\label{fig:volume-rendering-0p35}
	\end{subfigure}
	\begin{subfigure}{0.32\textwidth}
		\centering
		\includegraphics[trim={1cm 1 1cm 1 },clip,width=\textwidth]{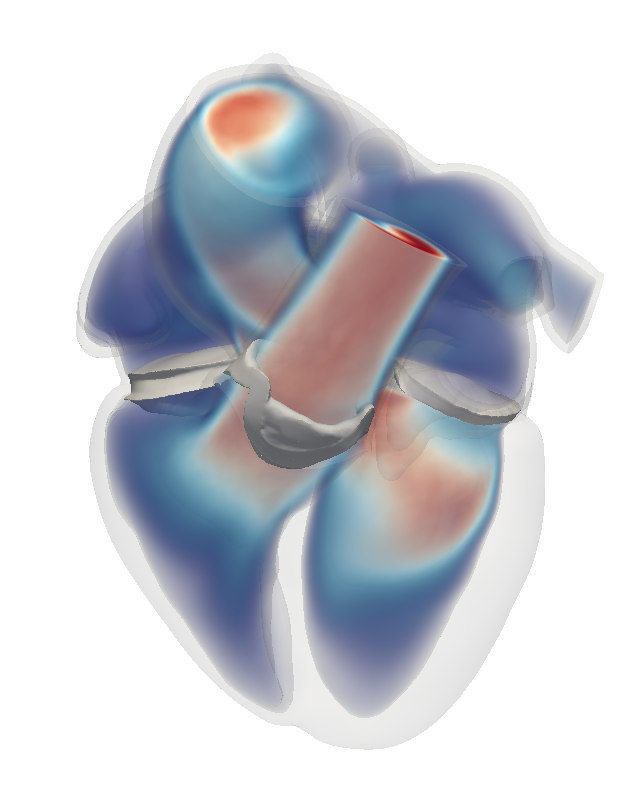}
		\caption{$t = 0.55$ s}
		\label{fig:volume-rendering-0p55}
	\end{subfigure}
	\begin{subfigure}{0.32\textwidth}
		\centering
		\includegraphics[trim={1cm 1 1cm 1 },clip,width=\textwidth]{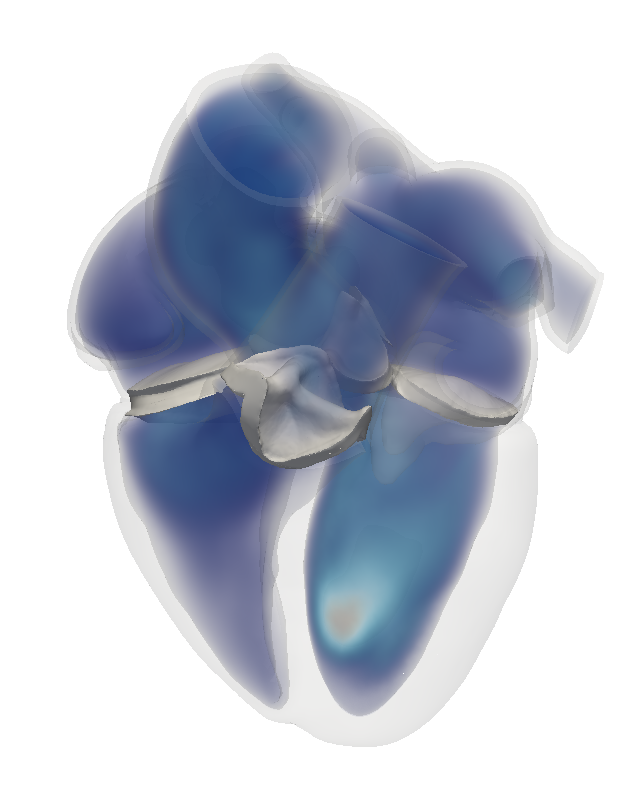}
		\caption{$t = 0.70$ s}
		\label{fig:volume-rendering-0p70}
	\end{subfigure}
	\begin{subfigure}{0.32\textwidth}
		\centering
		\includegraphics[trim={1cm 1 1cm 1 },clip,width=\textwidth]{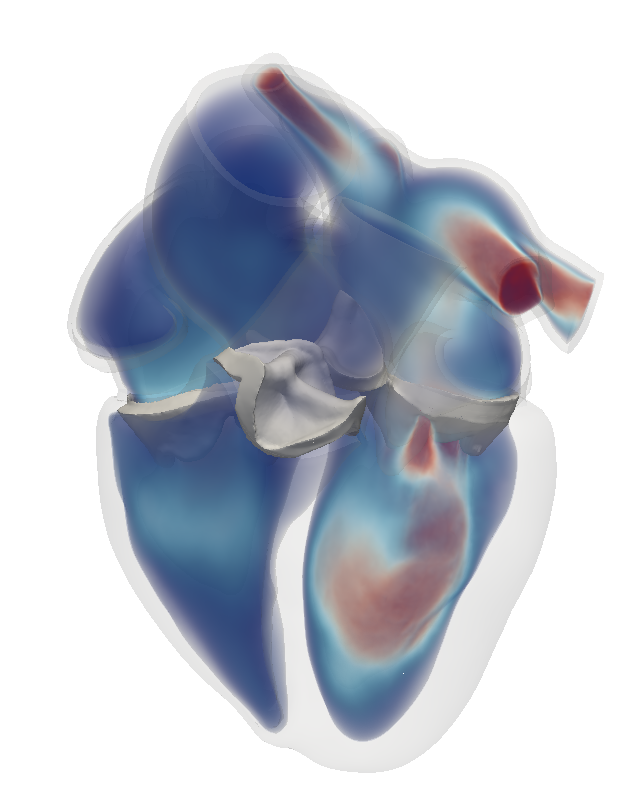}
		\caption{$t = 0.80$ s}
		\label{fig:volume-rendering-0p80}
	\end{subfigure}
	\\
	\centering
	\vspace{0.5cm}
	$|\bm u| $ [m/s]
	\\
	\begin{subfigure}{\textwidth}
		\centering
		\includegraphics[trim={1 1 1 1 },clip,width=0.8\textwidth]{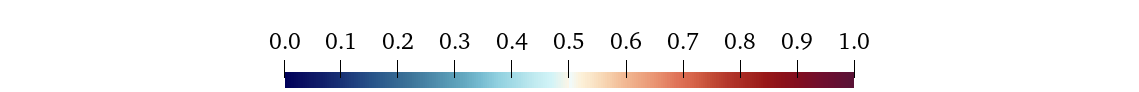}
	\end{subfigure}
	\caption{Volume rendering of velocity magnitude in different frames during the cardiac cycle: (a) diastolic a-wave peak, (b) isovolumetric contraction, (c) systolic peak, (d) mid systolic deceleration, (e) isovolumetric relaxation, (f) diastolic a-wave peak.}
	\label{fig:volume-rendering-velocity}
\end{figure}

During the ventricular passive filling, the blood flows from the pulmonary veins and the venae cavae {into} the LA and RA, respectively. Moreover, the atrioventricular valves are open and high-speed jets form between their leaflets. This moment corresponds to the E-wave of diastole (see \Cref{fig:velocities}). Consistently with clinical measurements, the computational model is able to correctly reproduce the formation of the clockwise jet in the LV, redirecting the blood towards the outflow tract~\cite{di2018jet, kilner2000asymmetric}. Furthermore, from \Cref{tab:validation-cfd}, we can observe that maximum velocity between MV and TV leaflets are in the physiological ranges. Furthermore, we also report average atrial pressure values and we found a general good agreement with reference data, even if left atrial pressure is slightly larger than our reference. \rev{To better appreciate how in silico results are in line with physiological ranges, we display each biomarker in \Cref{fig:cfd-normalized-graphic}. As for \Cref{fig:calibration-rdq20-graphic}, we normalize each biomarker in the interval $[-1, 1]$ with respect to the physiological values listed in \Cref{tab:validation-cfd}. }

\begin{table}
	\centering
	\begin{tabular}{c S c c c}
		\toprule
		\multicolumn{2}{c}{\textbf{Biomarker}} & \textbf{In silico} & \multicolumn{2}{c} {\textbf{Physiological values}} 
		\\
		\midrule
		Peak MV velocity & [\si{\metre\per\second}] & 1.03 &  0.89 $\pm$ 0.15 &~\cite{thomas1998peak}  
		\\
	    Peak AV velocity& [\si{\metre\per\second}] & 1.21 & 1.07 $\pm$ 0.18 &~\cite{bueno2006effect} 
		\\
		Peak TV velocity & [\si{\metre\per\second}] &  0.45 & 0.48 $\pm$ 0.11 &~\cite{choi2016normal} 
		\\
		Peak PV velocity & [\si{\metre\per\second}] &  1.15 & 0.80 to 1.20 &~\cite{ref_pulmonary_valve}  
		\\
		Mean LA pressure & [\si{\mmHg}] & 13.3 & 2 to 12 &~\cite{ref_pressures_chambers} 
		\\
		Peak LV pressure & [\si{\mmHg}] & 111 & 119 $\pm$ 13 &~\cite{sugimoto2017echocardiographic} 
		\\
		Mean RA pressure & [\si{\mmHg}] & 6.33 & 0 to 8 &~\cite{ref_pressures_chambers} 
		\\
		Peak RV pressure  & [\si{\mmHg}] & 37.0 & 35 $\pm$ 11 &~\cite{bishop1997clinical}
		\\
		\bottomrule
	\end{tabular}
	\caption{Fluid dynamics biomarkers obtained with the whole-heart CFD simulation (normal ranges or mean $\pm$ standard deviation). In silico values are computed by averaging fluid properties in spherical control volumes located between the valve leaflets (velocities) and in the heart chambers (pressures).}
	\label{tab:validation-cfd}
\end{table}

\begin{figure}
	\centering
	\begin{subfigure}{0.4\textwidth}
		\centering
		\includegraphics[trim={1 1 1 1 },clip,width=0.8\textwidth]{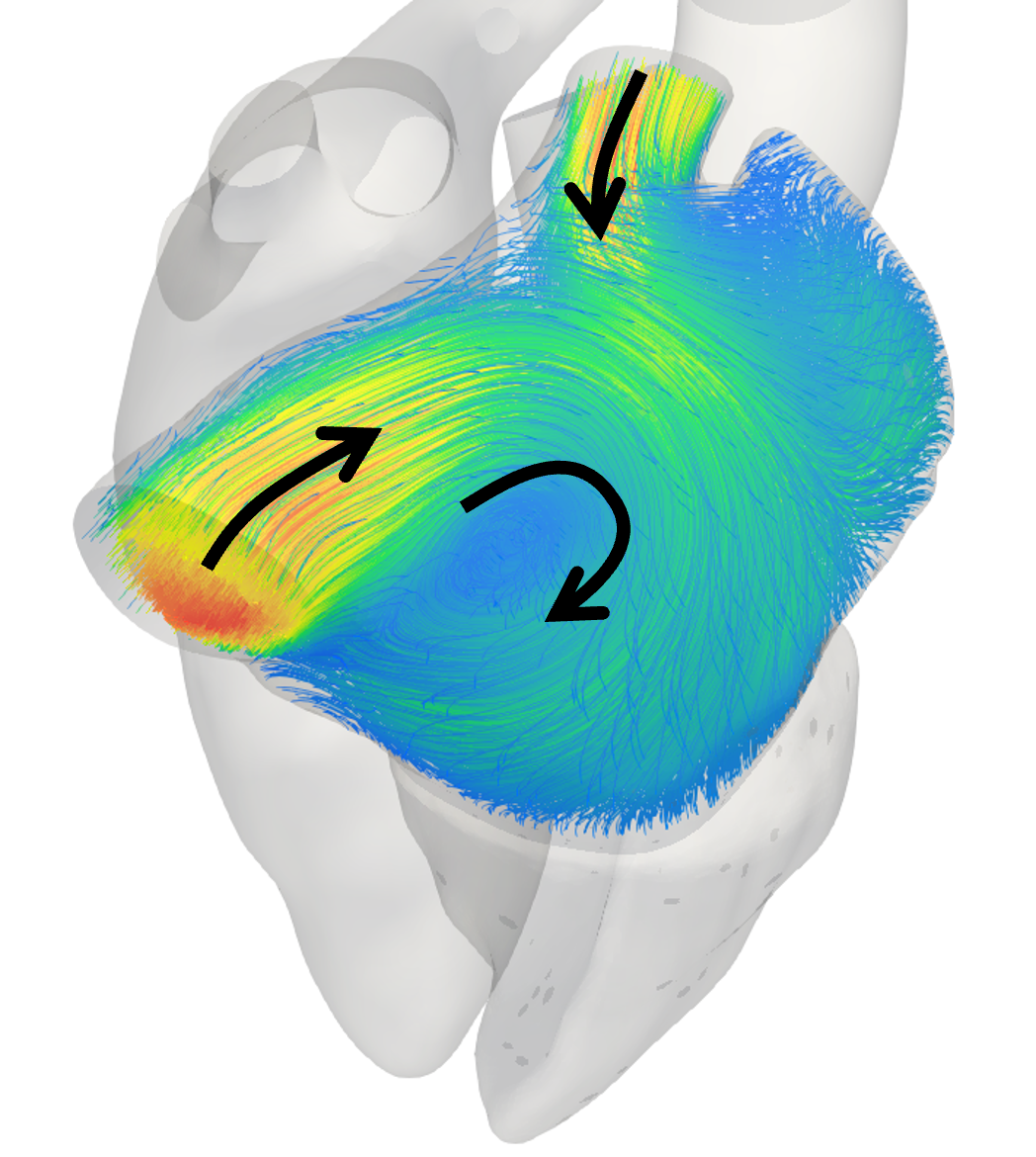}
		\caption{Right atrium $t=0.6$ s}
		\label{streamlines-ra-0p6}
	\end{subfigure}
	\begin{subfigure}{0.4\textwidth}
	\centering
	\includegraphics[trim={1 1 1 1 },clip,width=0.8\textwidth]{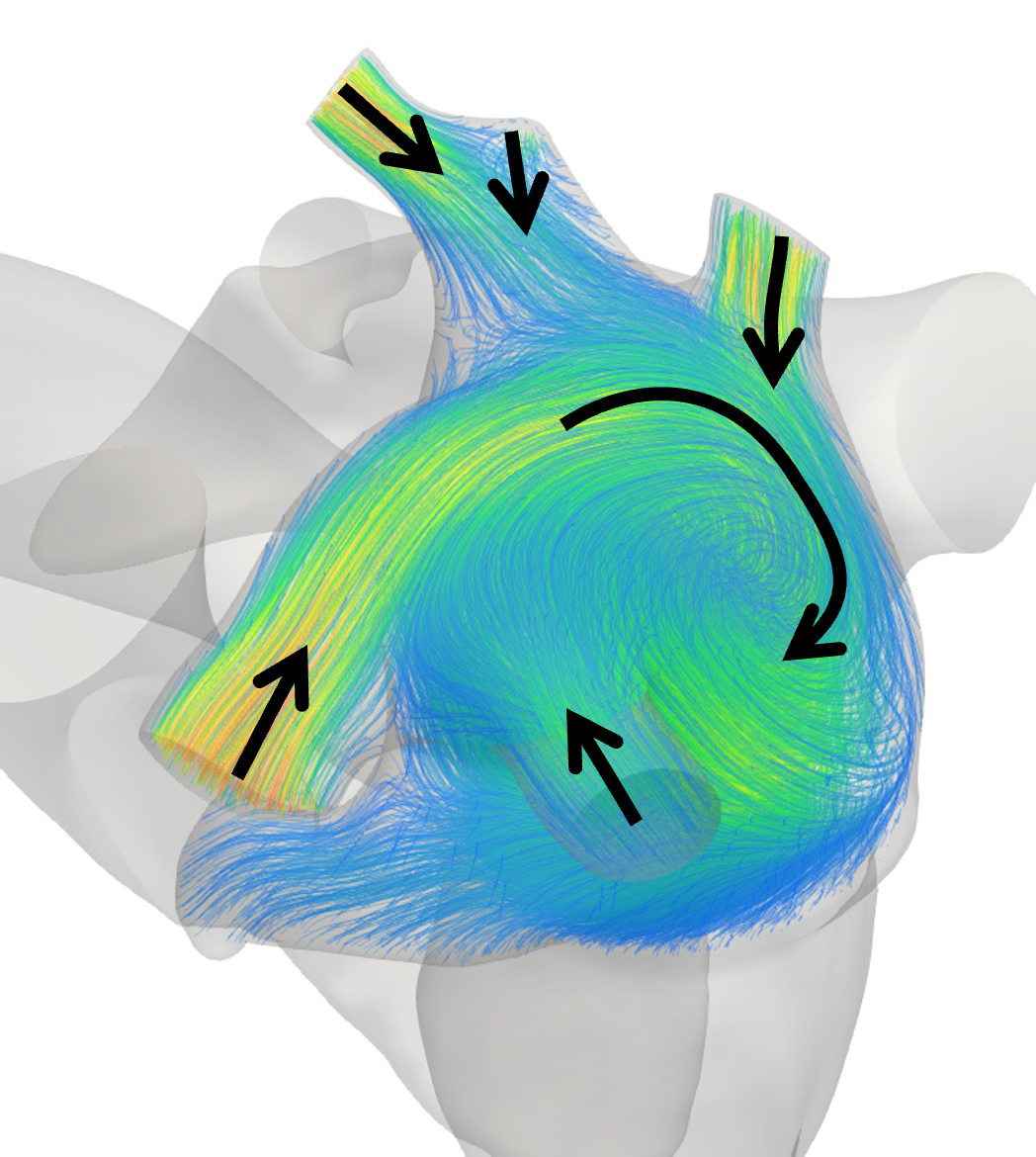}
			\caption{Left atrium $t=0.6$ s}
	\label{streamlines-la-0p6}
	\end{subfigure}
	\begin{center}
			\small{$|\bm u|$ [m/s]}
	\end{center}
	\begin{subfigure}{0.49\textwidth}
	\centering
	\includegraphics[trim={1 1 1 0.8cm},clip,width=0.8\textwidth]{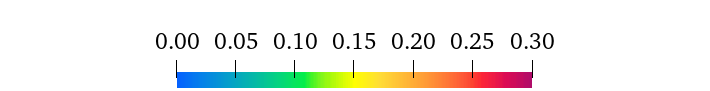}
	\end{subfigure}
	\\ \vspace{0.5cm}
	\begin{subfigure}{0.4\textwidth}
	\centering
	\includegraphics[trim={1 1 1 1 },clip,width=0.8\textwidth]{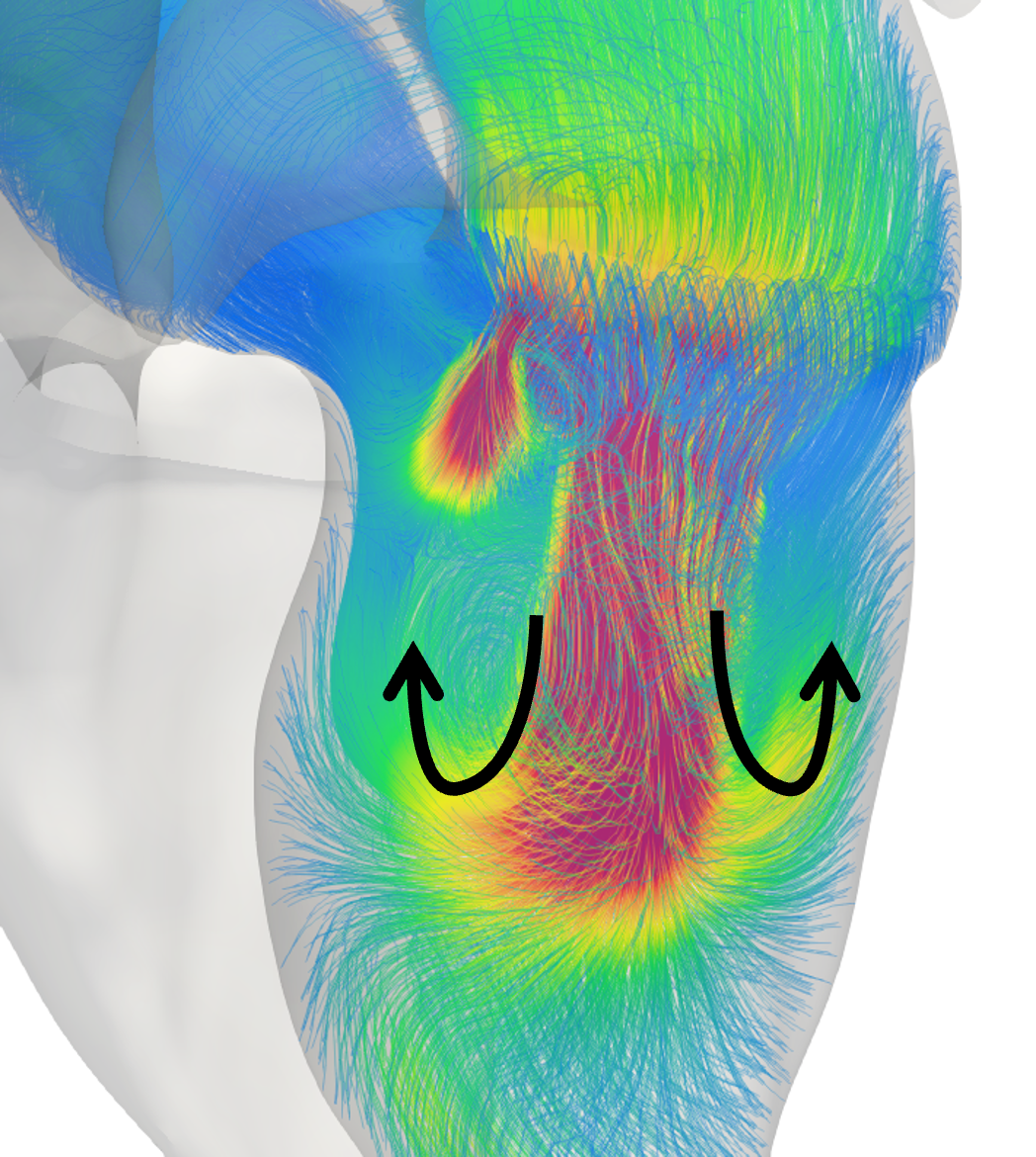}
			\caption{Left ventricle $t=0.775$ s}
	\label{streamlines-lv-0p775}
	\end{subfigure}
	\begin{subfigure}{0.4\textwidth}
	\centering
	\includegraphics[trim={1 1 1 1 },clip,width=0.8\textwidth]{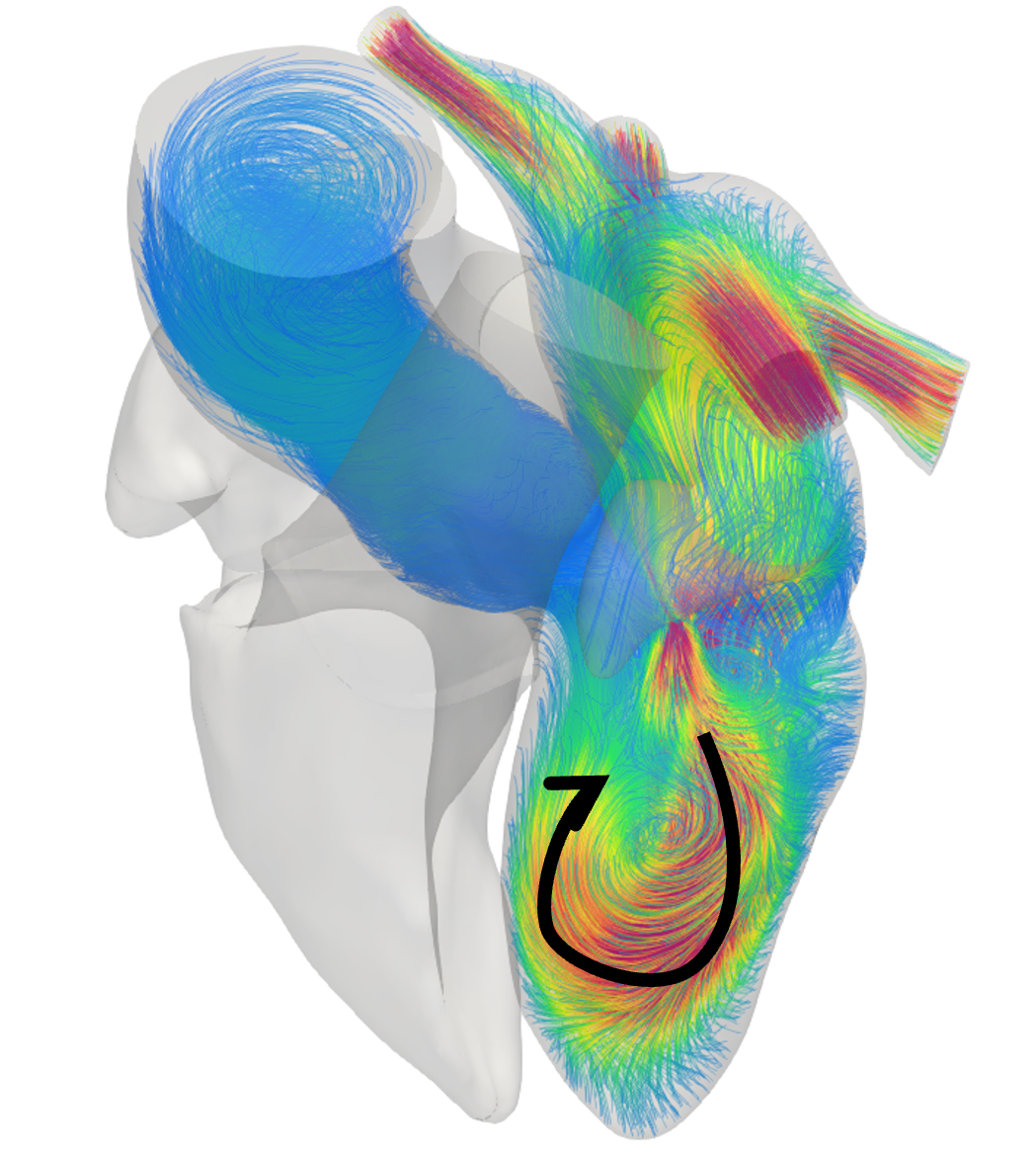}
			\caption{Left heart $t=0.8$ s}
	\label{streamlines-lh-0p8}
	\end{subfigure}
		\begin{center}
			\small{$|\bm u|$ [m/s]}
		\end{center}
	\begin{subfigure}{0.49\textwidth}
		\centering
		\includegraphics[trim={1 1 1 0.8cm},clip,width=\textwidth]{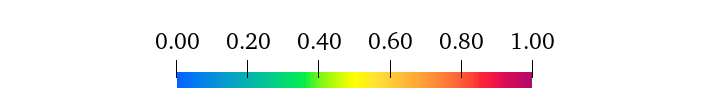}
	\end{subfigure}
	\caption{Streamlines colored according to velocity magnitude at different locations in the heart at different times: (a) right atrium during ventricular systole, (b) let atrium during ventricular systole, (c) clockwise and counter-clockwise vortices in the left ventricle during early diastole, (d) formation of clockwise jet in the left ventricle.}
	\label{fig:streamlines}
\end{figure}

\FloatBarrier
To better investigate blood flow patterns in the whole heart, we report the streamlines colored according to velocity magnitude in different chambers in \Cref{fig:streamlines}. We compare our in silico results with the MRI phase-velocity mapping visualization provided in Figure 1 of reference~\cite{kilner2000asymmetric}, and we found a good accordance. Specifically, \Cref{fig:streamlines}a shows the rotation of the blood in the RA as the chamber expands and the blood flows from the inferior and superior vena cava. {Similarly,} on the left side, the blood flows from the pulmonary veins to the expanding LA producing collision of blood jets and redirecting the flow towards the closed MV (see \Cref{fig:streamlines}b). During E-wave, as we show in \Cref{fig:streamlines}c, asymmetric recirculation is observed: shear layers roll through MV leaflets producing an O-vortex, as also seen in~\cite{chnafa2014image}. The counter-clockwise vortex under the posterior leaflet quickly disappears, and the clockwise vortex becomes larger and larger producing a clockwise jet, as described in~\cite{di2018jet}, and clearly observed in \Cref{fig:streamlines}d. 

\FloatBarrier

\subsection{{Application to a pathological scenario:} Left Bundle Branch Block {effects on hemodynamics}}
\label{sec:application}

\begin{figure}[t]
	\centering
	\begin{subfigure}{0.4\textwidth}
		\centering
		\includegraphics[trim={1 1 1 1 },clip,width=\textwidth]{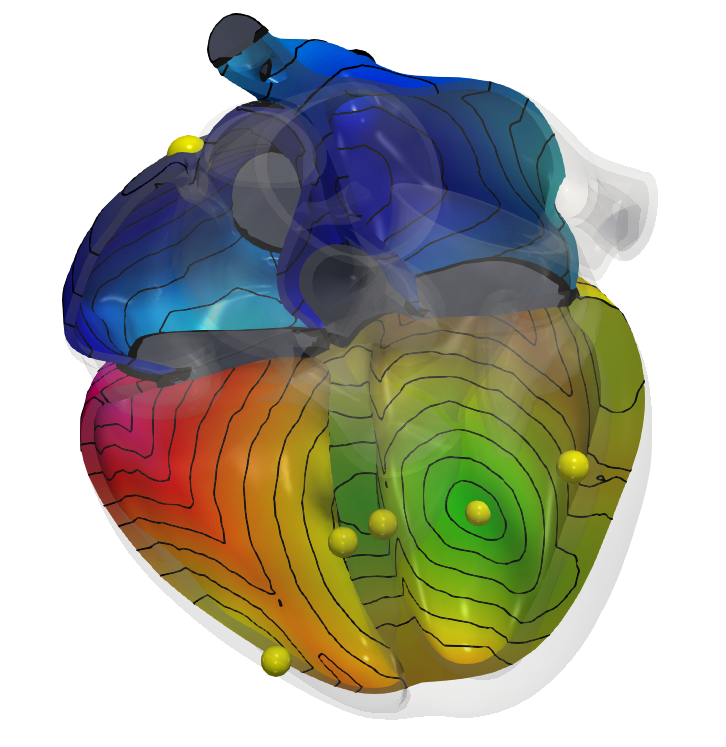}
		\caption{physiological}
		\label{fig:activation_physio}
	\end{subfigure}
	\begin{subfigure}{0.4\textwidth}
		\centering
		\includegraphics[trim={1 1 1 1 },clip,width=\textwidth]{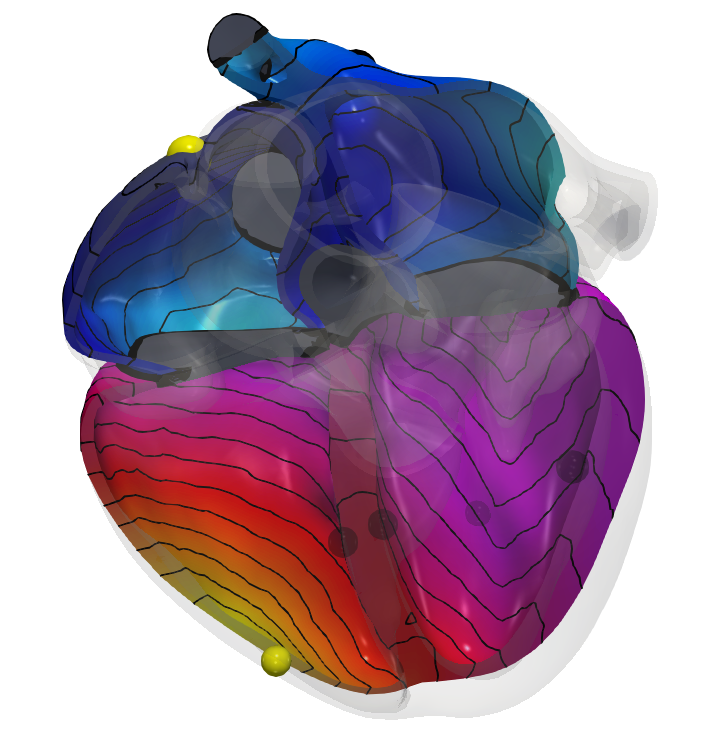}
		\caption{LBBB}
		\label{fig:activation_LBBB}
	\end{subfigure}
	\\
	\centering
	\vspace{0.5cm}
	\small{Activation time [ms]}
	\\
	\begin{subfigure}{0.8\textwidth}
		\centering
		\includegraphics[trim={1 1 1 0.5cm },clip,width=0.8\textwidth]{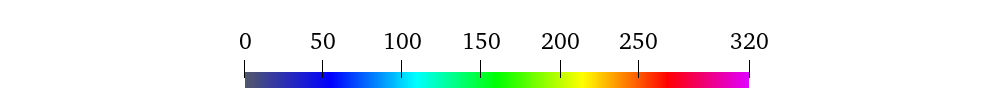}
	\end{subfigure}
	\caption{Activation maps and impulse sites (yellow spheres) in EM simulations; (a) physiological; (b) LBBB.}
	\label{fig:LBBB_activation}
\end{figure}

In this section, we apply our multiphysics computational model to investigate the hemodynamic consequences of the LBBB. This heart condition is commonly associated to an electrophysiological abnormality; however, it implies a cascade of adverse events due to the interaction among different {physical processes}.  LBBB consists in a slow or {even} absent conduction through the left bundle branch, causing a dyssynchronous contraction and relaxation of the left ventricle. Moreover, the LV dyssynchrony may have profound consequences on the heart hemodynamics, influencing flow patterns and, in turn, triggering heart remodeling \cite{littmann2000hemodynamic, eriksson2017left}. 

To simulate LBBB with our EM model, we deactivate the impulse sites in the LV and in the septum (see  \Cref{fig:LBBB_activation}), so that the signal is generated by the SAN and the RVm sites solely. Our modeling choice is consistent with the work of \cite{alessandrini2015pipeline}, where a severe case of LBBB is accounted for by activating only one site in the RV free wall.

\begin{figure}[t]
	\centering
	\includegraphics[trim={3cm 1 3cm 1 },clip,width=\textwidth]{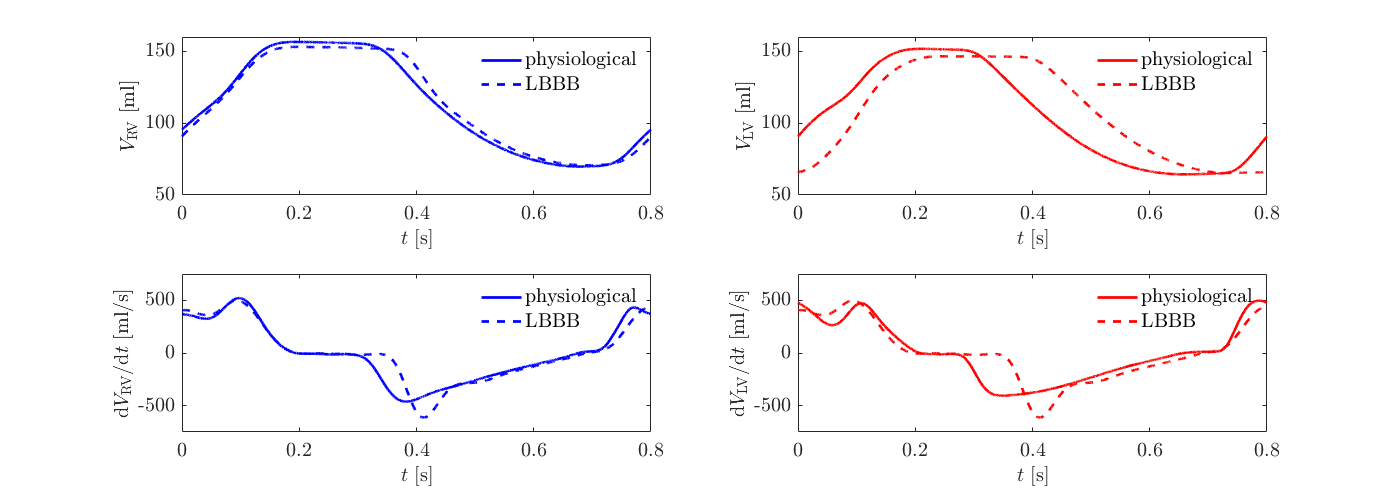}
	\caption{Ventricular volumes and ventricular volume derivatives for physiological and LBBB simulations.}
	\label{fig:volumes_volumesderivatives_physio_lbbb}
\end{figure}

\Cref{fig:volumes_volumesderivatives_physio_lbbb} {displays} volumes and their derivatives with respect to time {for} left and right ventricles, in {both} physiological and LBBB conditions. The electrical dyssynchrony between left and right parts produces a delay in the LV ejection and filling stages. Differently, no {significant} differences are observed in the RV volumes. Furthermore, compared to the physiological case (see \Cref{tab:RDQ20-results}) and consistently with \cite{littmann2000hemodynamic}, we measured reduced ejection fractions both in the right (\num{53.87}\%) and left (\num{55.39}\%) ventricles.
\begin{figure}[t]
    \centering
    \includegraphics[width=\textwidth]{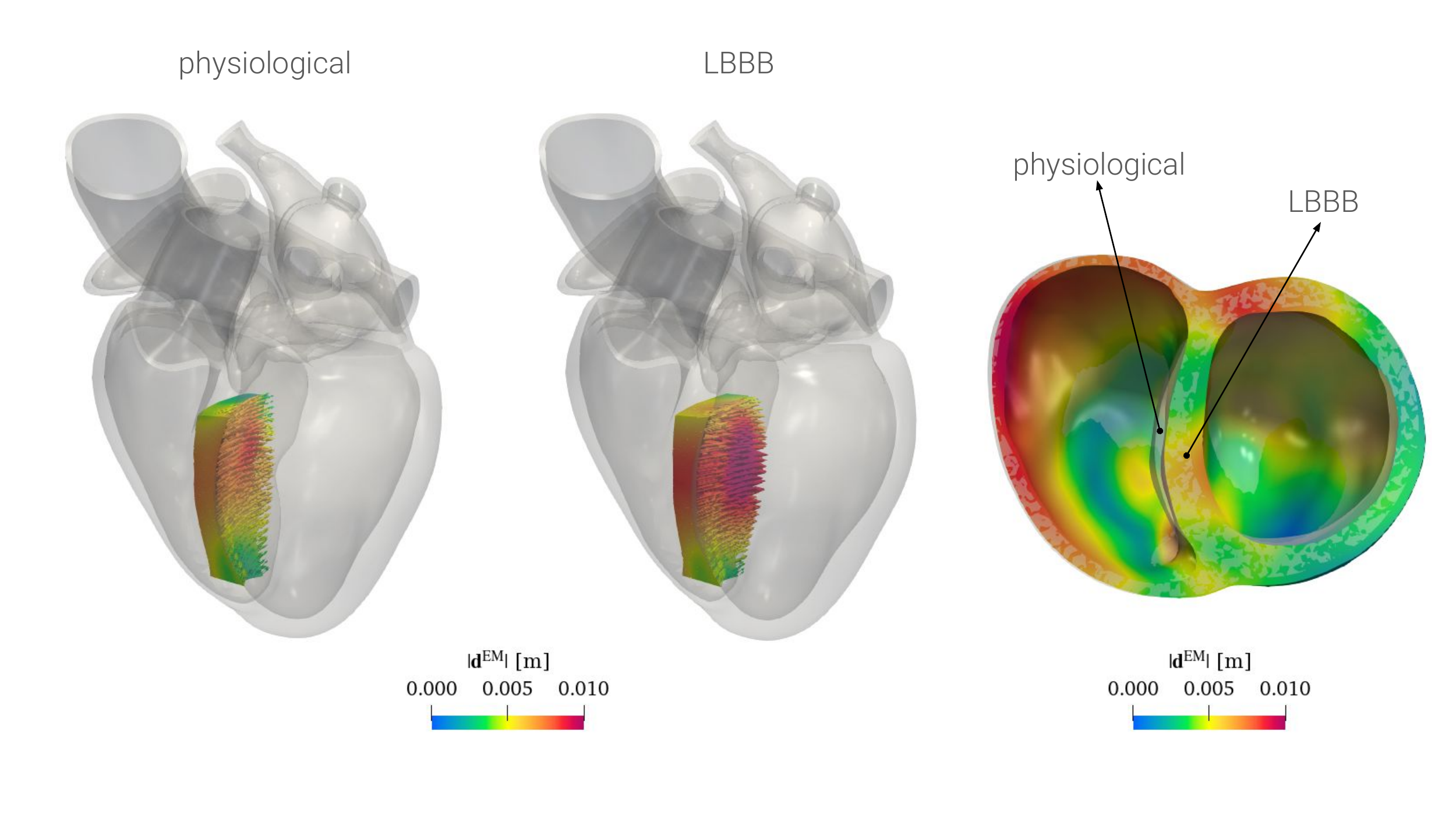}
    \caption{\rev{Comparison of the displacement field in the physiological and LBBB case at time $t=0.31$ s. From the left to the right: septum colored according to displacement magnitude with glyphs of displacement in the physiological and LBBB case; top-view of the left and right ventricles (in transparent the physiological case, colored according to displacement magnitude the LBBB case).}}
    \label{fig:septum-physio-vs-lbbb}
\end{figure}
\rev{Furthermore, we observe a different motion of the intraventricular septum during the pre-ejection phase. As we display in \Cref{fig:septum-physio-vs-lbbb}, the LBBB case is characterized by an abnormal leftward motion of the septum. Our result is consistent with different clinical observations \cite{dillon1974echocardiographic, gjesdal2011mechanisms}. }

\begin{figure}[t]
	\centering
	\begin{subfigure}{0.24\textwidth}
		\centering
		\includegraphics[trim={1 1 1 1 },clip,width=\textwidth]{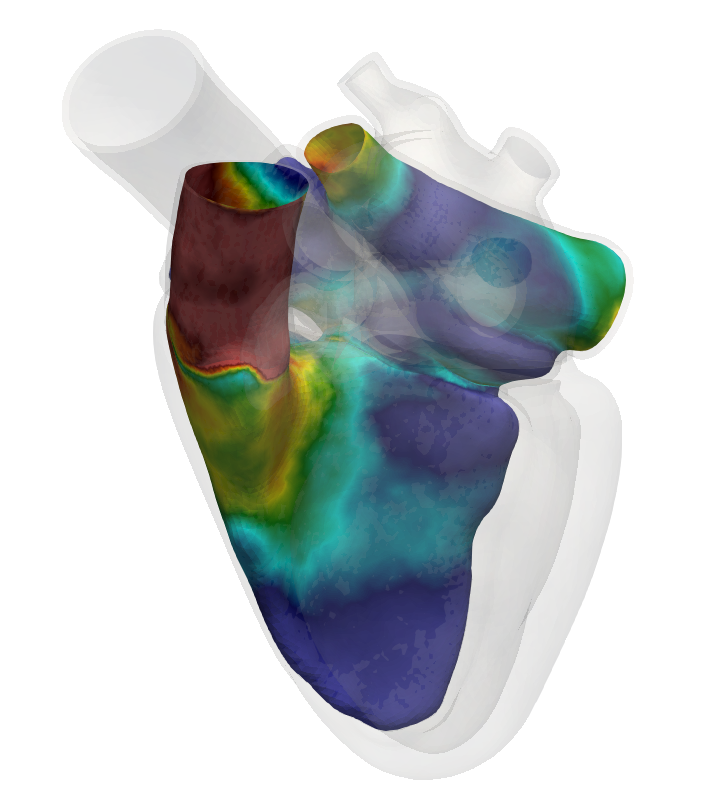}
		\caption{physiological, RH}
		\label{fig:TAWSS_physio_RH}
	\end{subfigure}
	\begin{subfigure}{0.24\textwidth}
		\centering
		\includegraphics[trim={1 1 1 1 },clip,width=\textwidth]{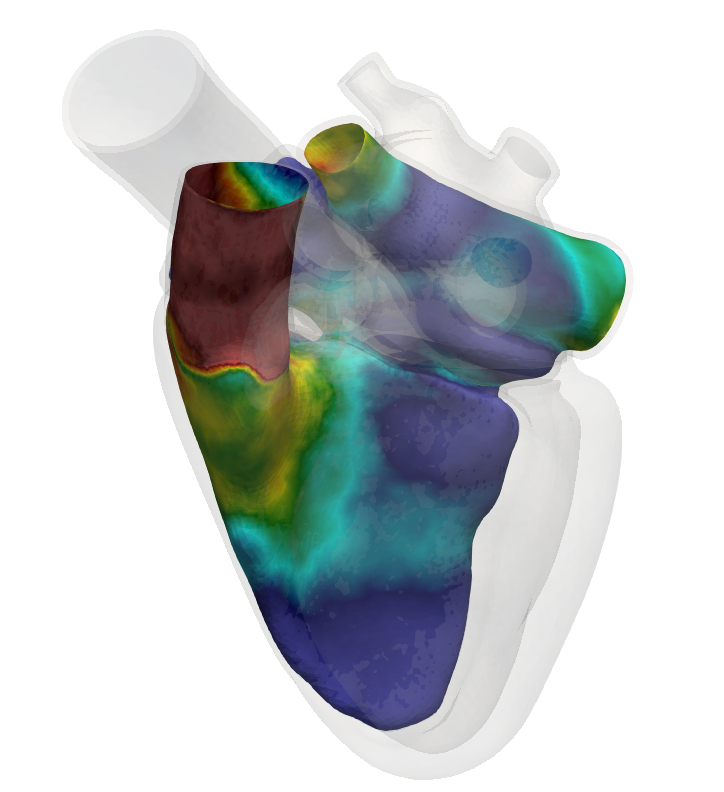}
		\caption{LBBB, RH}
		\label{fig:TAWSS_LBBB_RH}
	\end{subfigure}
	\begin{subfigure}{0.24\textwidth}
		\centering
		\includegraphics[trim={1 1 1 1 },clip,width=\textwidth]{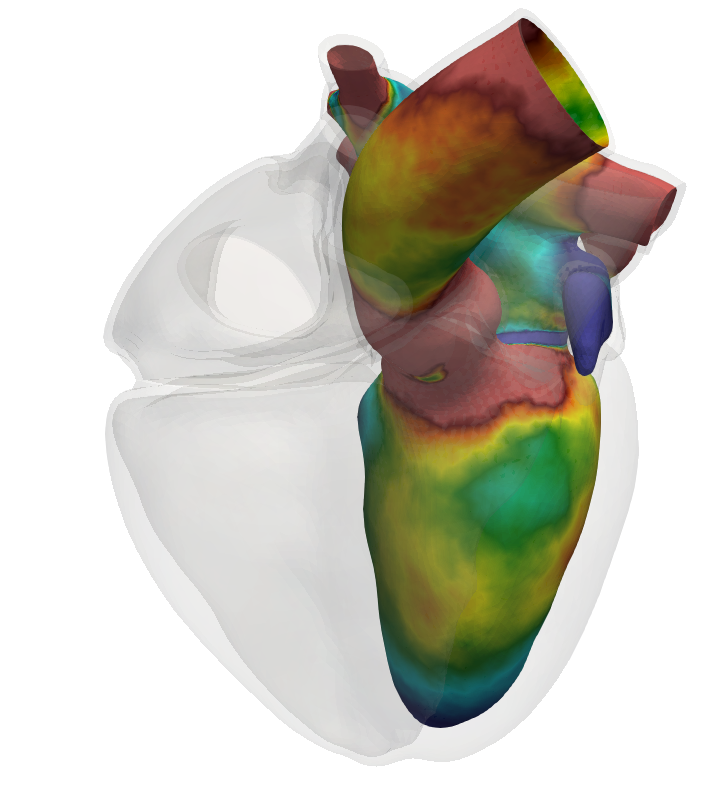}
		\caption{physiological, LH}
		\label{fig:TAWSS_physio_LH}
	\end{subfigure}
	\begin{subfigure}{0.24\textwidth}
		\centering
		\includegraphics[trim={1 1 1 1 },clip,width=\textwidth]{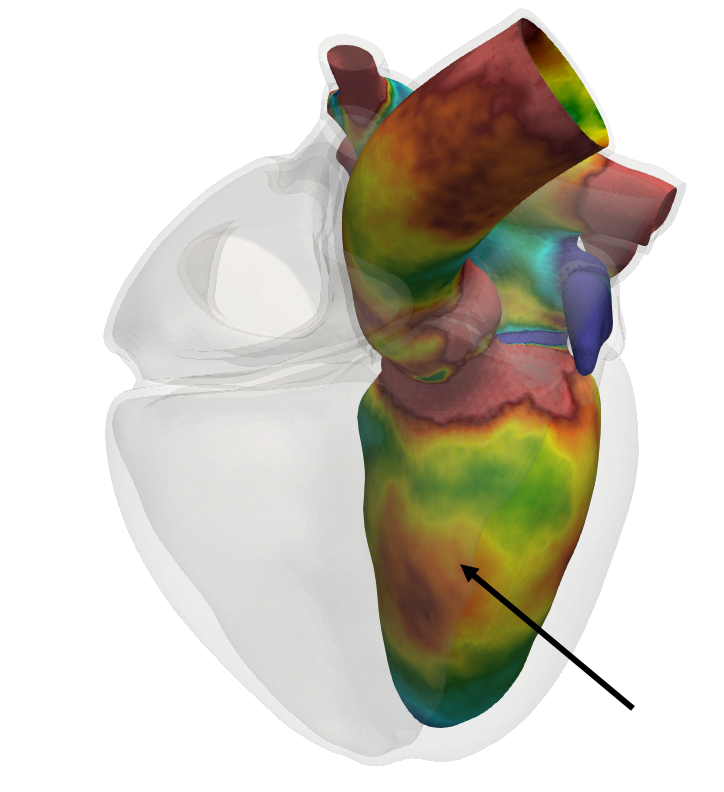}
		\caption{LBBB, LH}
		\label{fig:TAWSS_LBBB_LH}
	\end{subfigure}
	\\
	\centering
	\vspace{0.5cm}
	$\mathrm{TAWSS}$ [Pa]
	\\
	\begin{subfigure}{\textwidth}
		\centering
		\includegraphics[trim={1 1 1 1 },clip,width=\textwidth]{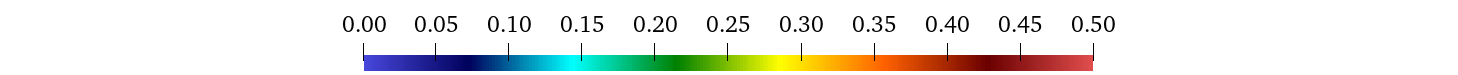}
	\end{subfigure}
	\caption{TAWSS of the whole heart in physiological and pathological conditions for the LH and RH. Results obtained with whole-heart EM-driven CFD simulations, the left and right parts are separated for visualization purposes. The largest differences are observed in the LV septum, as highlighted by the black arrow.}
	\label{fig:TAWSS}
\end{figure}
\begin{figure}[t]
	\centering
	\includegraphics[trim={1 3cm 1 1 },clip,width=\textwidth]{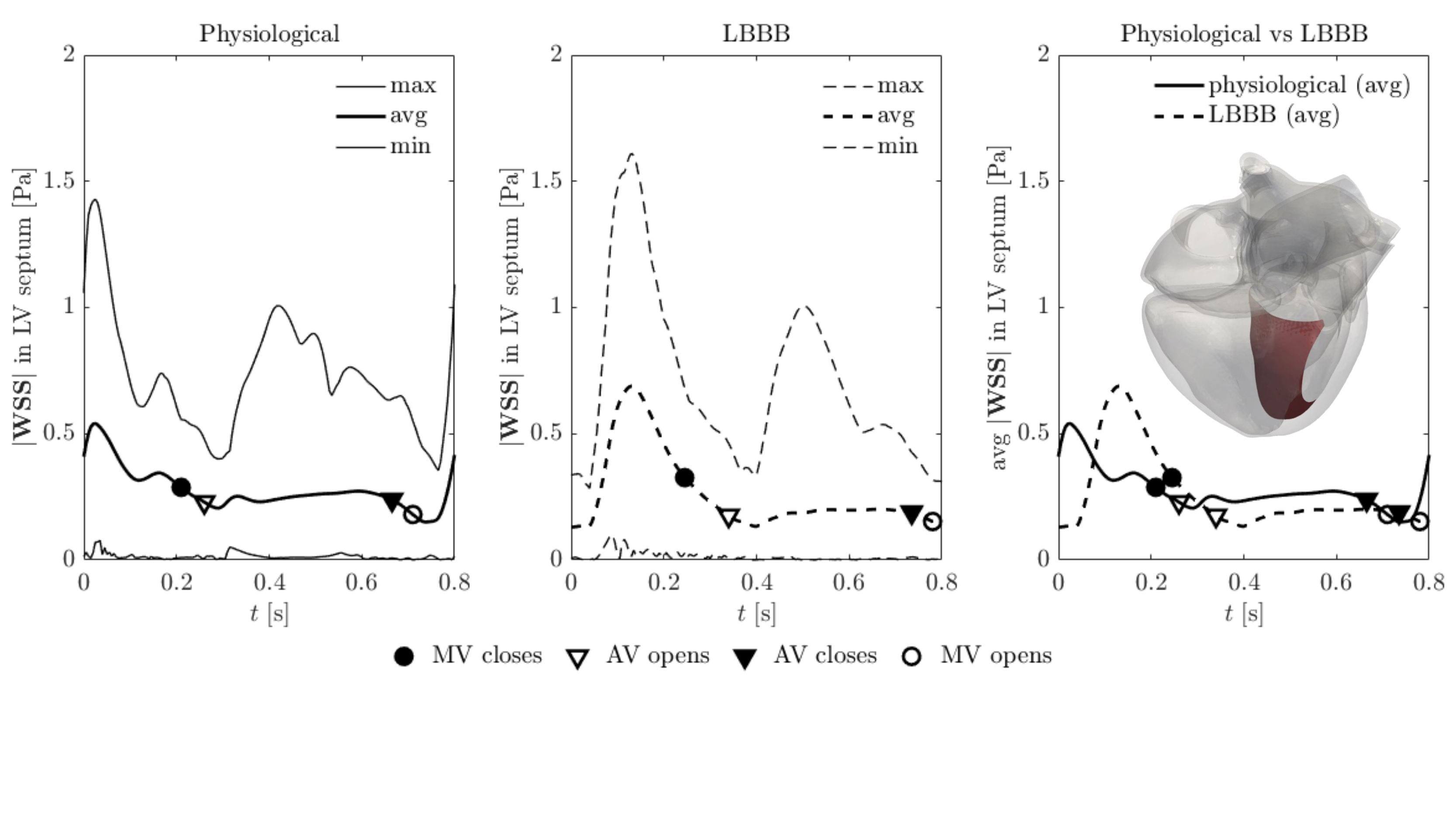}
	\caption{WSS in the LV septum over time. Maximum, average and minimum in physiological conditions (left); maximum, average and minimum in LBBB conditions (center); comparison between physiological and pathological conditions in terms of average values (right). The WSS is computed in the red portion shown on the right. {In terms of maximum values, the pathological peak is 12.7\% larger than the physiological case (cf. left and center figures). The average WSS peak increases of 27.4\% (right figure). }}
	\label{fig:LBBB_WSS_time}
\end{figure}

By means of our whole-heart EM driven CFD simulation, we can quantify how the pathology affects the endocardial wall stress. Let $\tau(\bm u) = 2\mu \epsilon(\bm u)$ be the viscous stress tensor, we compute the wall shear stress vector as
\begin{equation*}
	\mathbf{WSS}(\bm u) = \tau(\bm u) \bm n - (\tau(\bm u) \bm n\cdot \bm n)\bm n \qquad \text{ on } \partial \Omega_0 \times(0, T).
\end{equation*}
The time averaged wall shear stress (TAWSS) is then defined as
\begin{equation*}
	\mathrm{TAWSS} (\bm u) = \frac{1}{T} \int_{0}^T | \mathbf {WSS} (\bm u)|\, \mathrm d t  \qquad \text{ on } \partial \Omega_0.
\end{equation*}
\Cref{fig:TAWSS} shows the TAWSS in the right and left heart in physiological conditions and under LBBB. We notice that the TAWSS distribution is almost unchanged for the right heart. Conversely, we find that LBBB alters the wall shear stress in correspondence of the LV septum, suggesting the potential occurrence of remodeling phenomena. In \Cref{fig:LBBB_WSS_time}, we report the minimum, maximum and average WSS in the LV septum against time, for the physiological and LBBB simulations. During the ejection, the WSS values are similar, while significant differences are present during the filling phase. {As a matter of fact, under LBBB, the space-averaged WSS peak is 27.4\% higher than the physiological case (see \Cref{fig:LBBB_WSS_time}, right).} This is consistent with the work of Eriksson \textit{et al.} (2017), where they observed, by means of 4Dflow MRI data, that LV dyssynchronous motion influences blood flow patterns during diastole, contributing to the development of cardiac remodeling \cite{eriksson2017left}.

\FloatBarrier

\section{Discussion, limitations, and further developments}
\label{sec:limitations}
{We discuss some limitations of our study.} The computational model introduced cannot {fully represent} the isovolumetric phases in terms of pressure. Indeed, when both valves are closed and the ventricles are contracting/relaxing at constant volumes, the pressure is not well-defined, and thus prone to spurious oscillations. This is {due} to the fact that we are using kinematic conditions only in all the ventricular boundaries{. Indeed, we prescribe  the EM displacement on the endocardium and endothelium, and we model valves with the RIIS method using a penalty-based kinematic condition}. The {oscillatory} pressure during such phases does not influence the velocity field. However, it prevents from using the simulated pressure values to choose when to open and close the valves, forcing hence to prescribe a-priori opening and closing times. The use of bidirectionally coupled FSI models for the blood-myocardium or blood-valve systems (or for both of them) may allow to correctly capture the pressure transient during these phases, as seen for instance in \cite{bucelli2022mathematical} where a fully-coupled electro-mechano-fluid of the heart is considered. 

Moreover, we noticed that the ejection phase is too slow and the ventricular passive filling too fast, if compared with medical literature values. Consequently, E-wave and A-wave are characterized by comparable amplitudes, whereas the EA ratio should be approximately equal to 1.30 $\pm$ 0.570 \cite{thomas1998peak}. {In this respect, w}e believe that 
the use of ionic models with a more realistic decrease of calcium concentration is essential to better capture these phenomena.  

\st{Finally, we applied the computational model to a realistic, templated heart geometry. In order to move towards the realization of whole-heart digital twins, further developments should involve patient-specific cardiac simulations, accompanied with a stringent process of data assimilation, model validation and uncertainty quantification.}

\rev{A small time step is needed in the CFD simulation, particularly during the fast dynamics of the ejection phase. An adaptive time-stepping scheme would mitigate the computational burden of the overall multiphysics computational model. This scheme would involve using a small time step during ejection and larger values during the remaining phases of the cardiac cycle.}

\rev{Constructing a surrogate model would be a promising tool for expeditious calibration and sensitivity analysis. This model, incorporating inputs from both the EM and CFD models, along with outputs representing electrical, mechanical, and hemodynamical quantities of interest, would offer a mean to efficiently handle the complexity of the complex full order model. Previous works in the literature have successfully built such models for standalone cardiac electromechanics \cite{strocchi2023cell, regazzoni2022machine, longobardi2020predicting, salvador2023real} and standalone image-based hemodynamics \cite{karabelas2022global}.}

\rev{Ensuring the ``credibility evidence" for medical device submission of this intricate model would involve thorough verification and validation processes \cite{asme, fda_vvuq}. For simplified models, such as cardiac electrophysiology and hemodynamics, we presented verification studies in our recent releases of  \texttt{life$^\texttt{x}$} \cite{africa2023lifexep, africa2024lifexcfd}, our FE cardiac solver. Differently, for this multiphysics model, drawing inspiration from practices in cardiac electrophysiology \cite{niederer2011verification}, the creation of a simplified benchmark (e.g., an idealized left ventricle) incorporating the coupling of electrophysiology, mechanics, and hemodynamics would facilitate systematic verification of model's implementation.}

\rev{Finally, moving towards model validation and the development of whole-heart digital twins, future efforts should focus on patient-specific cardiac simulations. Integrating various in vivo data sources, such as electrocardiograms, cine cardiac magnetic resonance, and 4D flow MRI data, for comparison with in silico results presents a substantial challenge. However, overcoming obstacles in data assimilation, model parametrization, and validation would become imperative, given the inherent uncertainties associated with such in vivo data.}

\section{Conclusions}
\label{sec:conclusions}

In this paper, we introduced a computational model for the hemodynamics simulation of the whole human heart accounting for the main features affecting the intracardiac flows. We considered a realistic whole-heart geometry, and we employed a four-chamber 3D-0D electromechanical model to provide the displacement as input to the cardiac CFD model. We modelled the effect of cardiac valves in the fluid via a resistive immersed method and we accounted for transition-to-turbulence regime through the VMS-LES method. Moreover, for the first time, we coupled the 3D CFD model of  the whole heart to the surrounding closed-loop circulation, to get a geometric multiscale 3D-0D hemodynamic model of the entire cardiovascular system. We solved our multiphysics and multiscale computational model {using our} in-house finite element library \texttt{life}$^\texttt{x}$. 

We introduced a calibration of the \rev{active contraction model}, driving the electromechanical simulation, aimed at obtaining physiological realistic flowrates, and consequently blood velocities, in the CFD simulation. Our calibration highlights the effect of the parameters of the active force generation model, associated to the microscopic features of its kinematics and of the force-velocity relationship, on the macroscopic heartbeat indicators.

We carried out EM-driven CFD simulation on a realistic whole-heart geometry and we showed that the computational model can correctly reproduce blood velocities and pressure traces when we compare the results with clinical ranges from medical literature. Furthermore, we found that the computational model captures typical blood flow patterns observed in MRI phase-velocity mapping visualizations. 

Finally, we applied the whole-heart model to simulate the pathological scenario of Left Bundle Branch Block: we correctly predicted the electrical delay, the consequent mechanical dyssynchrony, a reduced ejection fraction, and an increasing wall shear stress in the left ventricular septum during the filling stage. Overall, this study confirms that {the} interaction of different physics in a high-fidelity integrated whole-heart model is essential {for} simulating the cardiac function, allowing to faithfully capture pathological events occurring at different physical levels.

\FloatBarrier

\appendix

\section{The open 0D circulation model}
\label{appendix:0dmodel}
The closed-loop lumped-parameter (0D) circulation model that we employ was proposed in \cite{regazzoni2022cardiac} and inspired by \cite{blanco20103d, hirschvogel2017monolithic}. To couple the 0D model to the 3D model of the heart, we follow the same steps presented in \cite{zingaro2022geometric}. Specifically, the open 0D system we get reads as follows: 
for any $t \in (0, T)$, 

\begin{subequations}
	\allowdisplaybreaks
	\begin{align}
		\frac{\mathrm dp_\text{AR}^\text{SYS}(t)}{\mathrm dt} & = \frac{1}{C_\text{AR}^\text{SYS}} \left (Q_\text{AV}(t) - Q_\text{AR}^\text{SYS}(t)
		\right ), \\
		\frac{\mathrm dp_\text{VEN}^\text{SYS}(t)}{\mathrm dt} & = \frac{1}{C_\text{VEN}^\text{SYS}} \left ( Q_\text{AR}^\text{SYS}(t) - Q_\text{VEN}^\text{SYS}(t) \right ),
		\\
		\frac{\mathrm dp_\text{AR}^\text{PUL}(t)}{\mathrm dt} & = \frac{1}{C_\text{AR}^\text{PUL}} \left (Q_\text{PV}(t) - Q_\text{AR}^\text{PUL}(t) \right ),
		\\
		\frac{\mathrm dp_\text{VEN}^\text{PUL}(t)}{\mathrm dt} & = \frac{1}{C_\text{VEN}^\text{PUL}} \left( Q_\text{AR}^\text{PUL}(t) - Q_\text{PUL}^\text{VEN, 3D}(t) \right ),
		\\
		\frac{\mathrm dQ_\text{AR}^\text{SYS}(t)}{\mathrm dt} & = \frac{R_\text{AR}^\text{SYS}}{L_\text{AR}^\text{SYS}} \left (- Q_\text{AR}^\text{SYS}(t) - \frac{p_\text{VEN}^\text{SYS}(t)-p_\text{AR}^\text{SYS}(t)}{R_\text{AR}^\text{SYS}} \right ),
		\\
		\frac{\mathrm dQ_\text{AR}^\text{PUL}(t)}{\mathrm dt} & = 	\frac{R_\text{AR}^\text{PUL}}{L_\text{AR}^\text{PUL}} \left (- Q_\text{AR}^\text{PUL}(t) - \frac{p_\text{VEN}^\text{PUL}(t)-p_\text{AR}^\text{PUL}(t)}{R_\text{AR}^\text{PUL}} \right ),
	\end{align}
	\label{eq:whole_heart_0D_3D_ODE}
\end{subequations}
solved with suitable initial conditions, and
\begin{subequations}
	\begin{align}
		p_\text{LA}^\text{in}(t) & = p_\text{VEN}^\text{PUL}(t) - R_\text{VEN}^\text{PUL}Q_\text{VEN}^\text{PUL}(t) - L_\text{VEN}^\text{PUL}\dfrac{\mathrm d Q_\text{VEN}^\text{PUL, 3D}(t)}{\mathrm dt}, \\
		p_\text{RA}^\text{in}(t) & = p_\text{VEN}^\text{SYS}(t) - R_\text{VEN}^\text{SYS}Q_\text{VEN}^\text{SYS}(t) - L_\text{VEN}^\text{SYS}\dfrac{\mathrm d Q_\text{VEN}^\text{SYS, 3D}(t)}{\mathrm dt}.
	\end{align}
	\label{eq:whole_heart_0D_3D_stationary}
\end{subequations}

\section{Setup of EM and CFD simulations}
\label{appendix:em-cfd}

We report in this section the values of the parameters used in the EM and CFD simulations.

\Cref{tab:circulation-em} reports the parameters for the circulation model used in the EM simulation. With respect to the model presented in \cite{fedele2022comprehensive}, we introduce two additional resistance elements ($R_\text{upstream}^\text{SYS}$ and $R_\text{upstream}^\text{PUL}$) and two inductance elements ($L^\text{AV}$, $L^\text{PV}$, see \Cref{fig:em-circulation-epicardium-fibers}c). They have the purpose of making the 0D circulation model, which surrogates the fluid dynamics, as similar as possible to the 3D CFD problem. For the same reason, the minimum resistances of the non-ideal diodes representing the AV and PV are increased with respect to \cite{fedele2022comprehensive}. 


The calibration of the RDQ20 \rev{active contraction model} is extensively discussed in \Cref{sec:calibration-rdq20mf}, and the corresponding parameter values are reported in \Cref{tab:RDQ20-setup}. As reference sarcomere length, we set $\mathrm{SL}_0$ = \num{2.2} \si{\micro\metre}. For all other parameters in the EM model, we use the same values as those of the baseline simulation of \cite{fedele2022comprehensive}.

\Cref{tab:circulation-cfd} reports the values of the initial circulation states for the CFD simulation. They are taken equal to the values reached at the beginning of the last heartbeat in the EM simulation. All resistance, capacitance and inductance parameters are the same as in the EM model (see \Cref{tab:circulation-em}).

\begin{table}
	\centering
	
	\begin{tabular}{l r c c}
		\toprule
		& \textbf{Parameter} & \multicolumn{2}{c}{\textbf{Value}} \\
		\midrule\multirow{6}{*}{Systemic arteries}
		& $R_\text{AR}^\text{SYS}$ & 0.48 & \si{\mmhg\second\per\milli\litre} \\
		& $C_\text{AR}^\text{SYS}$ & 1.50 & \si{\milli\litre\per\mmhg} \\
		& $L_\text{AR}^\text{SYS}$ & 0.005 & \si{\mmhg\square\second\per\milli\litre} \\
		& $R_\text{upstream}^\text{SYS}$ & 0.048 & \si{\mmhg\second\per\milli\litre} \\
		& $p_{\text{AR, 0}}^\text{SYS}$ & 83.9 & \si{\mmhg} \\
		& $Q_\text{AR, 0}^\text{SYS}$ & 0.0 & \si{\milli\litre\per\second} \\
		
		\midrule\multirow{5}{*}{Systemic veins}
		& $R_\text{VEN}^\text{SYS}$ & 0.26 & \si{\mmhg\second\per\milli\litre} \\
		& $C_\text{VEN}^\text{SYS}$ & 60 & \si{\milli\litre\per\mmhg} \\
		& $L_\text{VEN}^\text{SYS}$ & \num{5e-4} & \si{\mmhg\square\second\per\milli\litre} \\
		& $p_\text{VEN, 0}^\text{SYS}$ & 35.5 & \si{\mmhg} \\
		& $Q_\text{VEN, 0}^\text{SYS}$ & 0.0 & \si{\milli\litre\per\second} \\
		
		\midrule\multirow{6}{*}{Pulmonary arteries}
		& $R_\text{AR}^\text{PUL}$ & 0.032116 & \si{\mmhg\second\per\milli\litre} \\
		& $C_\text{AR}^\text{PUL}$ & 10 & \si{\milli\litre\per\mmhg} \\
		& $L_\text{AR}^\text{PUL}$ & 0.0005 & \si{\mmhg\square\second\per\milli\litre} \\
		& $R_\text{upstream}^\text{PUL}$ & 0.0032116 & \si{\mmhg\second\per\milli\litre} \\
		& $p_\text{AR, 0}^\text{PUL}$ & 14.90 & \si{\mmhg} \\
		& $Q_\text{AR, 0}^\text{PUL}$ & 0.0 & \si{\milli\litre\per\second} \\
		
		\midrule\multirow{5}{*}{Pulmonary veins}
		& $R_\text{VEN}^\text{PUL}$ & 0.035684 & \si{\mmhg\second\per\milli\litre} \\
		& $C_\text{VEN}^\text{PUL}$ & 16 & \si{\milli\litre\per\mmhg} \\
		& $L_\text{VEN}^\text{PUL}$ & 0.0005 & \si{\mmhg\square\second\per\milli\litre} \\
		& $p_\text{VEN, 0}^\text{PUL}$ & 13.58 & \si{\mmhg} \\
		& $Q_\text{VEN, 0}^\text{PUL}$ & 0.0 & \si{\milli\litre\per\second} \\
		
		\midrule\multirow{6}{*}{Valves}
		& $R_{\min}^\mathrm{MV}$ & 0.0075 & \si{\mmhg\second\per\milli\litre} \\
		& $R_{\min}^\mathrm{AV}$ & 0.0355 & \si{\mmhg\second\per\milli\litre} \\
		& $R_{\min}^\mathrm{TV}$ & 0.0075 & \si{\mmhg\second\per\milli\litre} \\
		& $R_{\min}^\mathrm{PV}$ & 0.0184 & \si{\mmhg\second\per\milli\litre} \\ 
		& $R_{\max}^\mathrm{MV}$, $R_{\max}^\mathrm{AV}$, $R_{\max}^\mathrm{TV}$, $R_{\max}^\mathrm{PV}$& 75006.2 & \si{\mmhg\second\per\milli\litre} \\
		& $L^\mathrm{AV}$, $L^\mathrm{PV}$ & \num{5e-4} & \si{\mmhg\square\second\per\milli\litre} \\
		
		\bottomrule
	\end{tabular}
	
	\caption{Parameters and initial conditions of the circulation model used for the EM simulation.}
	\label{tab:circulation-em}
\end{table}

\begin{table}
	\centering
	
	\begin{tabular}{l r c c}
		\toprule
		& \textbf{Parameter} & \multicolumn{2}{c}{\textbf{Value}} \\
		\midrule\multirow{2}{*}{Systemic arteries}
		& $p_\text{AR, 0}^\text{SYS}$ & 86.3480 & \si{\mmhg} \\
		& $Q_\text{AR, 0}^\text{SYS}$ & 109.6429 & \si{\milli\litre\per\second} \\
		
		\midrule\multirow{2}{*}{Systemic veins}
		& $p_\text{VEN, 0}^\text{SYS}$ & 34.4923 & \si{\mmhg} \\
		& $Q_\text{VEN, 0}^\text{SYS}$ & 112.9209 & \si{\milli\litre\per\second} \\
		
		\midrule\multirow{2}{*}{Pulmonary arteries}
		& $p_\text{AR, 0}^\text{PUL}$ & 22.2310 & \si{\mmhg} \\
		& $Q_\text{AR, 0}^\text{PUL}$ & 83.2132 & \si{\milli\litre\per\second} \\
		
		\midrule\multirow{2}{*}{Pulmonary veins}
		& $p_\text{VEN, 0}^\text{PUL}$ & 19.5813 & \si{\mmhg} \\
		& $Q_\text{VEN, 0}^\text{PUL}$ & 262.6397 & \si{\milli\litre\per\second} \\
		
		\bottomrule
	\end{tabular}
	
	\caption{Initial states of the circulation model for the CFD simulation. The values are equal to those of the circulation state variables at the beginning of the last heartbeat in the EM simulation. All remaining circulation parameters (resistances, capacitances, inductances) are the same as in the EM model (see \Cref{tab:circulation-em}).}
	\label{tab:circulation-cfd}
\end{table}

\clearpage

\section*{Acknowledgments}
AZ, LD and AQ received funding from the Italian Ministry of University and Research (MIUR) within the PRIN (Research projects of relevant national interest) 2017 “Modeling the heart across the scales: from cardiac cells to the whole organ” Grant Registration number 2017AXL54F). 

MB, RP, FR, LD and AQ acknowledge the ERC Advanced Grant iHEART, ``An Integrated Heart Model for the simulation of the cardiac function'', 2017–2023,  P.I. A. Quarteroni (ERC–2016– ADG, project ID: 740132). 
\begin{center}
	\includegraphics[width=0.7\textwidth]{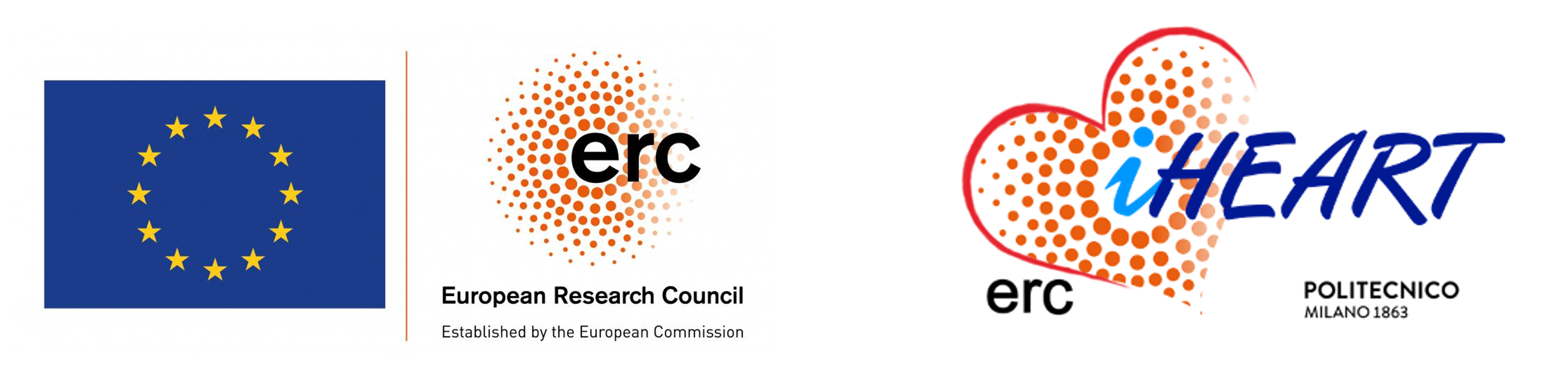}
\end{center}

The authors of this work are members of the INdAM group GNCS ``Gruppo Nazionale per il Calcolo Scientifico'' (National Group for Scientific Computing).

We acknowledge the CINECA award under the ISCRA initiative, for the availability of high performance computing resources and support under the projects IsC87\_MCH, P.I. A. Zingaro, 2021-2022 and \\ IsB25\_MathBeat, P.I. A. Quarteroni, 2021-2022. 

The present research is part of the activities of 650 ``Dipartimento di Eccellenza 2023-2027'', MIUR, Italy

\bibliographystyle{model1-num-names}
\bibliography{azingaro-jcp-bib}

\end{document}